\documentclass[a4paper, preprint, showkeys, nofootinbib]{revtex4}

\usepackage[ansinew]{inputenc}
\usepackage{amsmath, amsfonts} 
\usepackage{amsthm, amssymb}
\usepackage[english]{babel}
\usepackage{lmodern}
\usepackage{graphicx}
\usepackage{subfigure}

\theoremstyle{definition}

\newtheorem{remark}{Remark}

\usepackage{verbatim}
\usepackage{array}
\usepackage{float}
\usepackage{subfigure}
\usepackage{color}

\newcommand{\guillemets}[1]{``#1''}
\newcommand{\ve}[1]{\mathbf{#1}}
\newcommand{\red}[1]{\textcolor{black}{#1}}

\begin{document}

\title{Extreme phase sensitivity in systems with fractal isochrons}
\author{A. Mauroy}
\email{a.mauroy@ulg.ac.be}
\affiliation{Department of Electrical Engineering and Computer Science, University of Liège, 4000 Liège, Belgium}
\author{I. Mezi\'c}
\email{mezic@engr.ucsb.edu}
\affiliation{Department of Mechanical Engineering, University of California Santa Barbara, Santa Barbara, CA 93106, USA}

\begin{abstract}
Sensitivity to initial conditions is usually associated with chaotic dynamics and strange attractors. However, even systems with (quasi)periodic dynamics can exhibit it.
In this context we report on the fractal properties of the isochrons of some continuous-time asymptotically periodic systems.
We define a global measure of phase sensitivity that we call the {\it phase sensitivity coefficient} and show that it is an invariant of the system related to the capacity dimension of the isochrons.
Similar results are also obtained with discrete-time systems. As an illustration of the framework, we compute the phase sensitivity coefficient for popular models of bursting neurons, \red{suggesting that some} elliptic bursting neurons are characterized by isochrons of high fractal dimensions and exhibit a very sensitive (unreliable) phase response.
\end{abstract}


\maketitle

\section{Introduction}

Isochrons and asymptotic phase play a central role for the study of asymptotically periodic systems. The isochrons have been introduced in \cite{Winfree2} as the sets of initial states that converge to the same trajectory on the limit cycle. Equivalently, they are the sets of states that share the same \emph{asymptotic phase} \cite{Winfree}. The notions of isochrons and asymptotic phase are of paramount importance to capture the system sensitivity to external perturbations \cite{Taylor,Sacre3}. Indeed, an external perturbation has a significant permanent effect on the system only if the perturbed and unperturbed trajectories lie on different isochrons (associated with a noticeable asymptotic phase difference). Beyond their applications to sensitivity analysis, isochrons and asymptotic phase also lead to a powerful phase reduction of the dynamics that is widely used for studying synchronization properties of coupled limit-cycle systems (see e.g. \cite{Brown,Ermentrout_book,Hoppensteadt,Kuramoto_book}).

Thanks to a numerical method based on the so-called Koopman operator framework and introduced in \cite{Mauroy_Mezic}, it has been observed recently that isochrons may exhibit a fractal geometry \cite{MRMM_bursting}. (Note that a similar phenomenon was discussed in \cite{MezicBana} for the asymptotic phase of a discrete-time map.) \red{The fractal property of the isochrons is naturally explained in presence of transient chaos \cite{Yorke_chaotic_transient}, a regime characterized by sensitivity to initial conditions caused by a nonattracting (chaotic) set \cite{chaotic_dyn_book}. Also, this property is related to an extremely high phase sensitivity of the system, which has a dramatic effect on the response to external inputs and perturbations.}

\red{The main goal of this paper is to propose a theoretical framework---and associated numerical methods---that complements the brief and empirical observations given in \cite{MRMM_bursting} on the fractal properties of the isochrons.} To that end, we define \emph{the phase sensitivity coefficient} that quantifies the overall \red{uncertainty of the asymptotic phase under small perturbations}. Using a result of \cite{McDonald}, we prove that this coefficient is closely related to the capacity dimension of the isochrons \cite{Farmer_fractal_dim}. A consequence of our results is that, when the isochrons are fractal, a significant decrease of the intensity of a noise perturbation can only slightly reduce the average uncertainty on the phase.

\red{In contrast to local measures of sensitivity (e.g. finite-time Lyapunov exponents), the phase sensitivity coefficient defined in the present paper is an invariant of the system that captures a global property. An important motivation of our approach is therefore} to provide a framework for comparing the overall phase sensitivity of different asymptotically periodic systems. As an illustration, the framework is applied to several types of bursting neurons. The results show that elliptic bursting neurons \red{tested} are highly sensitive (i.e. with fractal isochrons of high dimension) and therefore characterized by unreliable responses to external inputs. Their (finite) phase response curve is also shown to be fractal.

The paper is organized as follows. In Section \ref{prelimin}, we rigorously define the notions of asymptotic phase and isochrons. Section \ref{sec_first_obs} presents basic observations and descriptions of the fractal properties of the isochrons. The relationship between phase sensitivity and fractal dimension of the isochrons is discussed in Section \ref{sec_sensitivity_fractal_dim}. We apply the results to bursting neuron models in Section \ref{sec_bursting}. Finally, concluding remarks are given in Section \ref{conclu}.

\section{Asymptotic phase and isochrons}
\label{prelimin}

In this section, we introduce the concepts of phase function and isochrons, for both continuous-time and discrete-time systems. A numerical method is also presented for the computation of the phase function.

\subsection{Continuous time}

We consider a nonlinear system $\dot{\ve{x}}=\ve{F(x)}$, with $\ve{x} \in \mathbb{R}^N$ and $\ve{F}$ analytic, which generates a flow $\varphi:\mathbb{R}^+ \times \mathbb{R}^N \rightarrow \mathbb{R}^N$ and admits a periodic orbit $\Gamma$ of period $T_0=2\pi/\omega_0$, i.e. $\varphi(T_0,\ve{x}^\gamma)=\ve{x}^\gamma$ with $\ve{x}^\gamma\in \Gamma$. Each point $\ve{x}^\gamma$ of the periodic orbit is associated with a phase $\theta \in \mathbb{S}^1$ according to the mapping $\ve{x}^\gamma(\theta)=\varphi((\theta/2\pi)T_0,\ve{x}^\gamma_0)$, where $\ve{x}^\gamma_0=\ve{x}^\gamma(0)$ is an arbitrarily chosen point of $\Gamma$ \cite{Winfree}.

If $\Gamma$ is an asymptotically stable limit cycle with a basin of attraction $\mathcal{B}\subseteq \mathbb{R}^N$, the (asymptotic) phase function $\Theta:\mathcal{B} \rightarrow \mathbb{S}^1$ assigns the same phase \red{$\theta \in \mathbb{S}^1$} to the initial states converging to the same trajectory on the limit cycle, i.e.
\begin{equation}
\label{phase_fct}
\Theta(\ve{x})=\theta \quad \Leftrightarrow \quad \lim_{t\rightarrow \infty} \|\varphi(t,\ve{x})-\varphi(t,\ve{x}^\gamma(\theta))\|=0\,.
\end{equation}
The level sets of $\Theta$---i.e. the sets of states that share the same asymptotic behavior---are the so-called isochrons \cite{Winfree2}
\begin{equation*}
\mathcal{I}_\theta=\{\ve{x} \in \mathcal{B} | \Theta(\ve{x})=\theta\} \,.
\end{equation*}
If the limit cycle is normally hyperbolic, the isochrons $\mathcal{I}_\theta$ are co-dimension-$1$ manifolds that invariantly foliate the basin of attraction \cite{Hirsch}. Figure \ref{vdp_iso}(a) shows the asymptotic phase and $5$ isochrons for the Van der Pol model. 

\red{Since the phase function is defined only for initial conditions of trajectories converging to the limit cycle, it is not defined outside the basin of attraction $\mathcal{B}$. In particular, the boundary of $\mathcal{B}$ is called the \emph{phaseless set} $\mathcal{S}$ and satisfies the following property \cite{Guckenheimer_iso}: the values of the phase function evaluated on any neighborhood of $\mathcal{S}$ span the entire circle $\mathbb{S}^1$ (equivalently, the isochrons come arbitrarily close to $\mathcal{S}$). The phaseless set usually corresponds to a (unstable) fixed point and its stable manifold. In Figure \ref{vdp_iso}(a), the phaseless set is the unstable fixed point at the origin.}

Primarily developed in the context of computational neuroscience, an important tool using the notion of phase is the well-known (finite) phase response curve
\begin{equation}
\label{PRC}
Z_{\ve{e}}(\theta)=\Theta(\ve{x}^\gamma(\theta)+\ve{e})-\theta\,,
\end{equation}
with $\ve{x}^\gamma(\theta) \in \Gamma$ and $\ve{e} \in \mathbb{R}^N$. It represents the phase shift of a trajectory on the limit cycle that is subjected to an impulsive perturbation $\ve{u}(t)=\ve{e} \delta(t)$, where $\delta(t)$ is the Dirac function.

\subsection{Discrete time}

The phase function $\Theta$ can also be defined in the case of discrete-time systems. Consider the map $\ve{x}(t+1)=\ve{F}(\ve{x}(t))$, $\ve{x} \in \mathbb{R}^N$, $t\in \mathbb{N}$, which generates the flow $\varphi:\mathbb{N} \times \mathbb{R}^N \rightarrow \mathbb{R}^N$. We assume that the system admits an invariant set $\Gamma$, a topological circle, on which the dynamics has an irrational rotation number $\nu_0$.
This invariant set can be characterized as a closure of a dense orbit, i.e. $\Gamma=\overline{\cup_{t\in \mathbb{N}} \, \varphi(t,\ve{x})}$ for some $\ve{x}\in \mathbb{R}^N$. In this case, there exists a sequence $\{t_k\}_{1\leq k \leq \infty}$ such that $\lim_{k \rightarrow \infty} t_k \omega_0 \bmod 2\pi = 0$, with $\omega_0=2 \pi \nu_0$, and $\lim_{k \rightarrow \infty} \varphi(t_k,\ve{x}^\gamma)=\ve{x^\gamma}$ for all $\ve{x}^\gamma\in \Gamma$.

\begin{remark}[Proof of the existence of the sequence $t_k$] We have by definition
\begin{equation*}
\omega_0 = \lim_{t\rightarrow \infty}\frac{F^t(s(\ve{x}^\gamma))-s(\ve{x}^\gamma)}{t}\,,
\end{equation*}
where $s:\Gamma \rightarrow [0,2\pi]$ defines a coordinate on the circle $\mathbb{S}^1([0,2\pi])$ and $F$ is the lifting of the map on the circle to real line.
Let the sequence ${\epsilon_k}$ converge to $0$ and define
\begin{equation}
\label{nu_k}
\omega_k=\frac{F^{t_k}(s(\ve{x}^\gamma))-s(\ve{x}^\gamma)}{t_k},
\end{equation}
where $F^{t_k}(s(\ve{x}^\gamma))-s(\ve{x}^\gamma) \bmod 2\pi<\epsilon_k$ or $F^{t_k}(s(\ve{x}^\gamma))-s(\ve{x}^\gamma)\bmod 2\pi>2 \pi-\epsilon_k$ (we know such $t_k$ exist due to the fact that every trajectory in the circle is dense). Therefore, it follows from \eqref{nu_k} that $(\omega_k t_k) \bmod 2\pi<\epsilon_k$ or $(\omega_k t_k) \bmod 2\pi>2\pi-\epsilon_k$ and thus $\lim_{k\rightarrow \infty}(\omega_0 t_k) \bmod 2\pi=0$.
\end{remark}

\begin{figure}[h]
\begin{center}
\subfigure[]{\includegraphics[width=7cm]{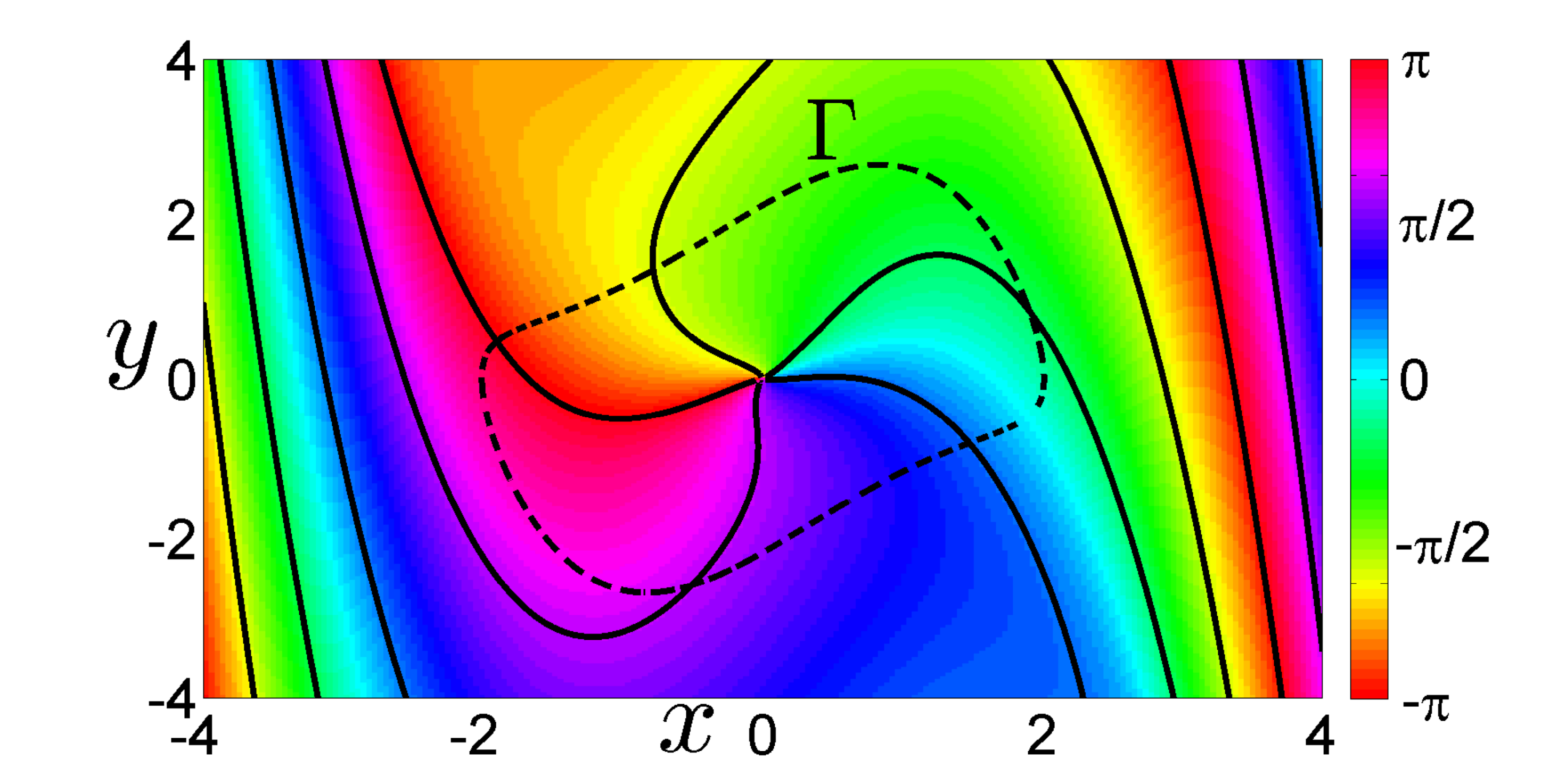}}
\subfigure[]{\includegraphics[width=7cm]{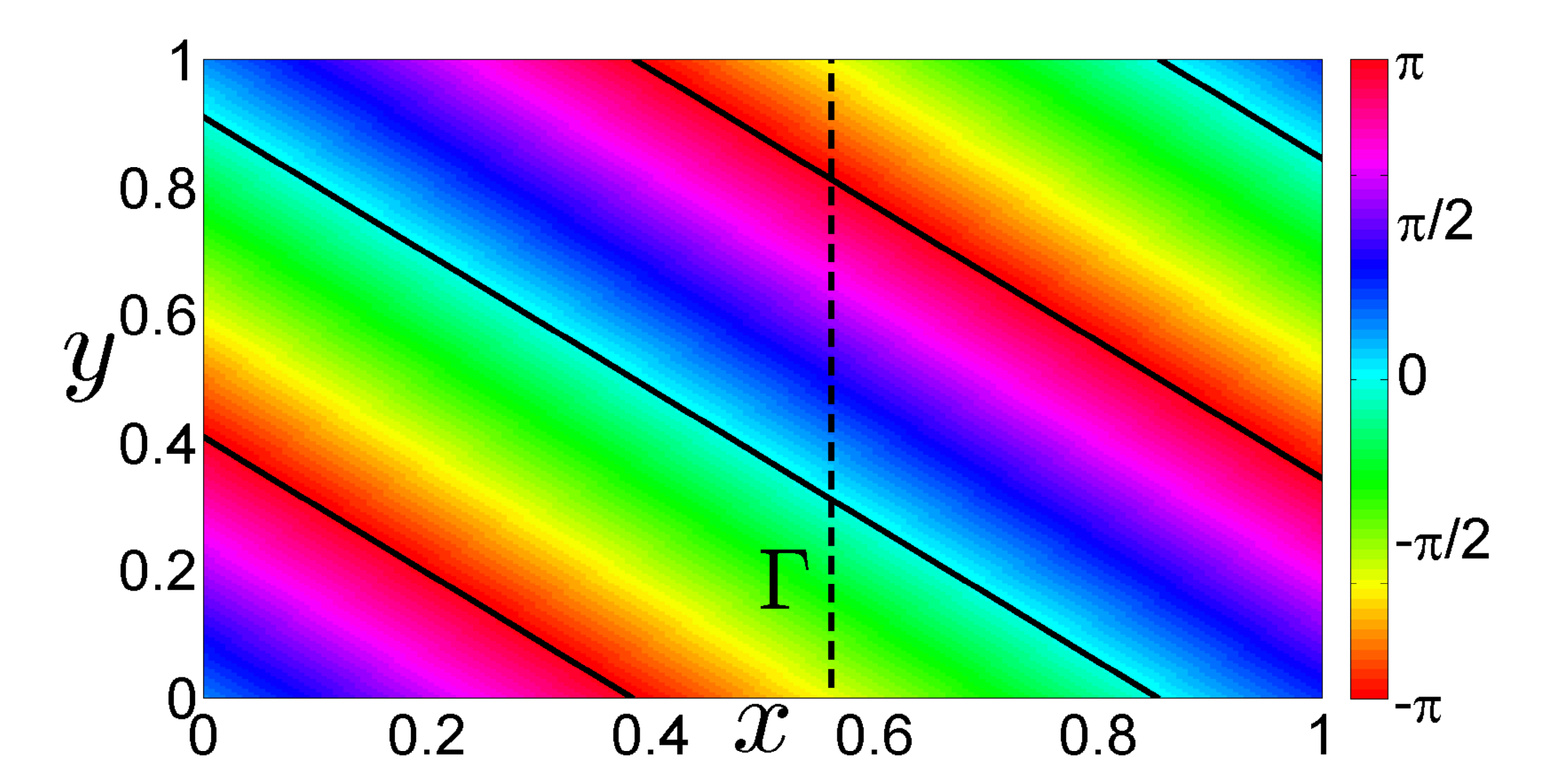}}
\caption{(a) Asymptotic phase and $5$ isochrons (equally spaced in phase) for the Van der Pol model $\dot{x}=y$, $\dot{y}=(1-x^2)y-x$. The unstable fixed point at the origin is the phaseless set $\mathcal{S}$. \red{The black dashed line is the limit cycle $\Gamma$.} (b) Asymptotic phase and $2$ isochrons for the discrete-time map (taken from \cite{MezicBana}) $x(t+1)=\gamma(x(t)-\nu_0)+\nu_0 \bmod 1$, $y(t+1)=y(t)+x(t) \bmod 1$ with $\nu_0=0.5613245623$, $\gamma=0.06123456756432$. \red{The black dashed line is the attractor $\Gamma$.}}
\label{vdp_iso}
\end{center}
\end{figure}

The phase is defined as follows. Each point $\ve{x}^\gamma \in \Gamma$ is associated with a phase $\theta \in \mathbb{S}^1([0,2\pi])$ according to the mapping $\ve{x}^\gamma(\theta)=\lim_{k \rightarrow \infty} \varphi(t_k,\ve{x}^\gamma_0)$, where $\ve{x}^\gamma_0=\ve{x}^\gamma(0)$ is an arbitrarily
chosen point of $\Gamma$ and where the sequence $\{t_k\}_{1\leq k \leq \infty}$ satisfies $\lim_{k \rightarrow \infty} t_k \omega_0 \bmod 2\pi = \theta$. If the invariant set $\Gamma$ is \red{an attractor with a basin of attraction $\mathcal{B}$ (i.e. $\Gamma$ is the smallest set such that $\overline{\cap_{t\in \mathbb{N}} \, \varphi(t,\ve{x})} \subseteq \Gamma$ for all $\ve{x}\in \mathcal{B}$)}, then the phase function is defined on $\mathcal{B}$ by \eqref{phase_fct}. In addition, the isochrons can be defined as the level sets of the phase function, but they might not be connected manifolds. The phase function and 2 isochrons of a simple discrete-time map are shown in Figure \ref{vdp_iso}(b).

\subsection{Numerical computation}

There exist several methods for computing the isochrons of limit cycles (see e.g. \cite{Guillamon,Huguet,Izi_book,Langfield,Osinga,Sherwood}). For our purpose, we will use the forward-integration method proposed recently in \cite{Mauroy_Mezic}, which is based on the fact that the phase function is related to an eigenfunction of the so-called Koopman operator \cite{MezicBana}. (Note that the method presented in \cite{Guillamon,Huguet} precisely solves a partial differential equation which is closely related to the eigenvalue equation for the Koopman operator.) The method is appropriate to deal with the complex geometry of the isochrons that we investigate in this paper. It is also well-suited to the computation of isochrons in non-planar models. Through this framework, the phase function of continuous-time systems is directly given by the argument $\angle$ of the Fourier average evaluated along the trajectories, i.e.
\begin{equation}
\label{Fourier_av}
\Theta(\ve{x})=\angle \left(\lim_{T \rightarrow \infty} \frac{1}{T} \int_0^T g \circ \varphi(t,\ve{x}) \, e^{-i\omega_0 t} \, dt \right) \,,
\end{equation}
where $g \in C^1:\mathcal{B} \rightarrow \mathbb{C}$ is an arbitrary function (observable) such that the first Fourier coefficient (first harmonic) of the periodic function $g \circ \varphi(t,\ve{x}^\gamma)$ (with $\ve{x}^\gamma \in \Gamma$) is nonzero. Note that the state $\ve{x}^\gamma_0 \in \Gamma$ associated with the phase $\theta=0$ is determined by the specific choice of $g$. Similarly, the phase function of a discrete-time map is given by 
\begin{equation}
\label{Fourier_av2}
\Theta(\ve{x})=\angle \left(\lim_{T \rightarrow \infty} \frac{1}{T} \sum_{t=0}^T g \circ \varphi(t,\ve{x}) \, e^{-i \omega_0 t} \, \right) \,.
\end{equation}

The Fourier averages \eqref{Fourier_av} and \eqref{Fourier_av2} can be easily computed through the numerical integration of trajectories, with the initial conditions on a uniform grid that spans a region of interest in the state space. The isochrons are obtained by plotting the level sets of these Fourier averages. Further numerical details can be found in Appendix \ref{app_numerical}

\section{Fractal properties of the phase function}
\label{sec_first_obs}

We consider particular (continuous-time and discrete-time) dynamical systems and show that their associated phase function exhibits fractal patterns, suggesting that the isochrons are fractal. These systems are characterized by a very high phase sensitivity.

\subsection{Continuous-time model}
\label{sec_continuous_model}

The  Lorenz system
\begin{equation}
\label{Lorenz}
\begin{array}{rcl}
\dot{x} & = & \sigma (y-x) \\
\dot{y} & = & x(r -z) - y  \\
\dot{z} & = & xy - b z
\end{array}
\end{equation}
admits a stable limit cycle for the parameters $\sigma=10$, $b=8/3$, and $r=320$. \red{It also exhibits transient chaos \cite{Yorke_chaotic_transient} and has a chaotic saddle (see e.g. \cite{chaotic_dyn_book}, also called fractal repeller \cite{Gaspard_fractal_repeller}), i.e. a fractal invariant set containing the union of all unstable periodic orbits. The chaotic saddle and its stable manifold correspond to the phaseless set $\mathcal{S}$, since trajectories starting on this set do not converge toward the limit cycle. The chaotic saddle is fractal, and therefore $\mathcal{S}$ is also characterized by fractal properties. In addition, $\mathcal{S}$ includes the stable manifold of the saddle node at the origin, which was shown in \cite{McDonald} to have a fractal Cantor-like geometry.}

\begin{figure}[h]
\begin{center}
\subfigure[~Phase function in the cross-section $y=50$]{\includegraphics[width=7cm]{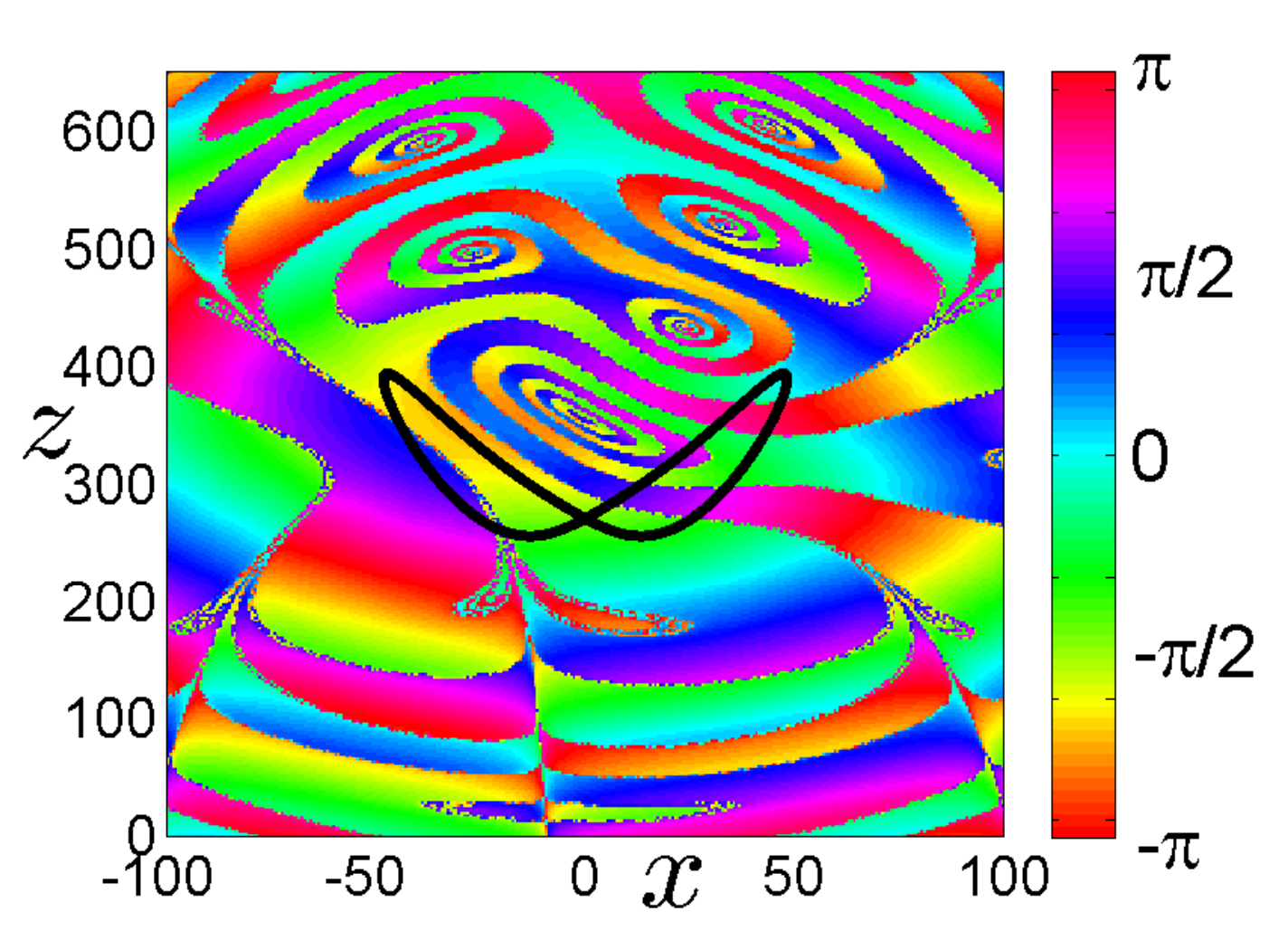}}
\subfigure[~Phase function in the cross-section $x=25$]{\includegraphics[width=7cm]{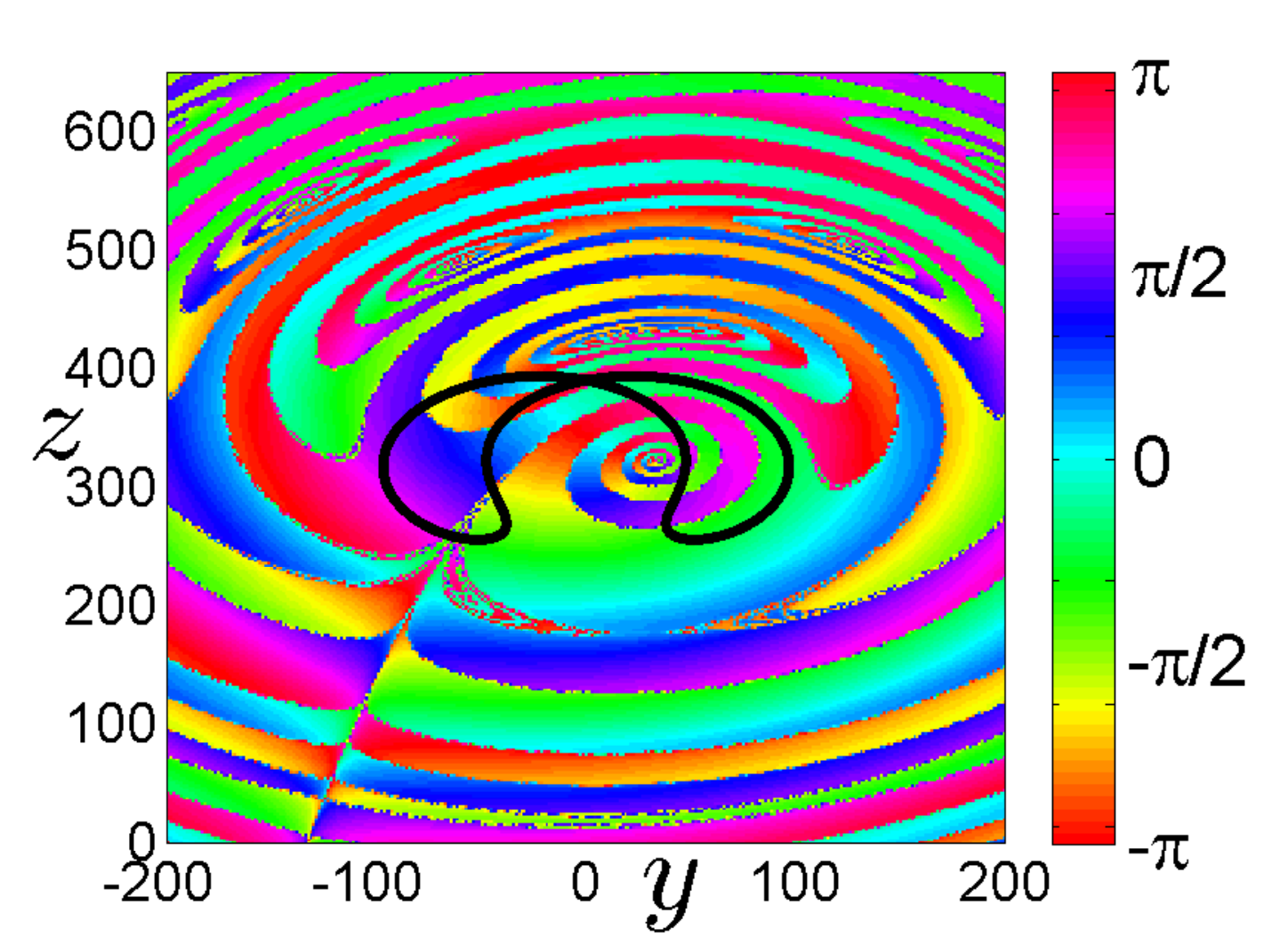}}
\subfigure[~Phase function in the cross-section $z=319$]{\includegraphics[width=7cm]{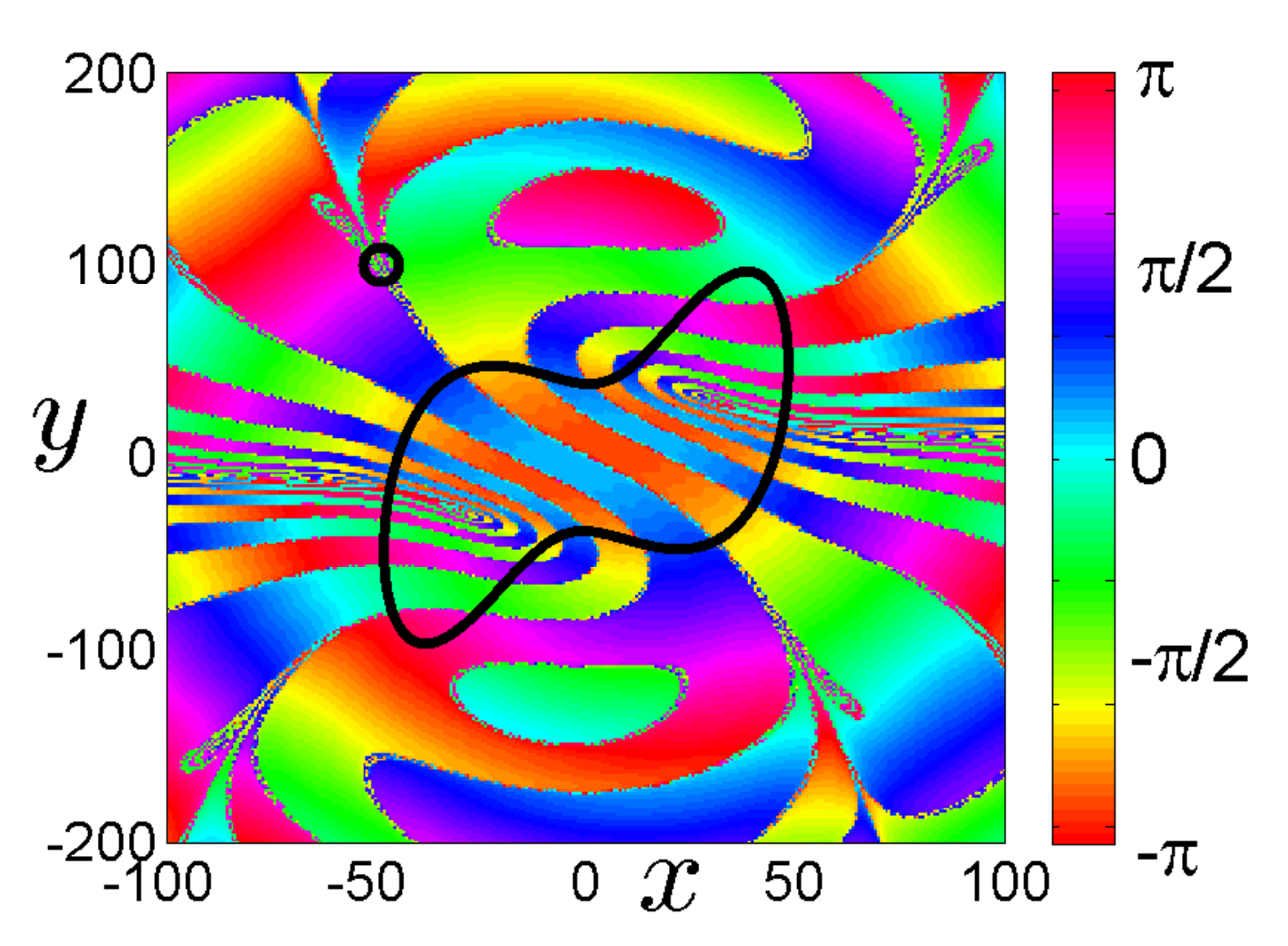}}
\subfigure[~Close-up in the cross-section $z=319$]{\includegraphics[width=7cm]{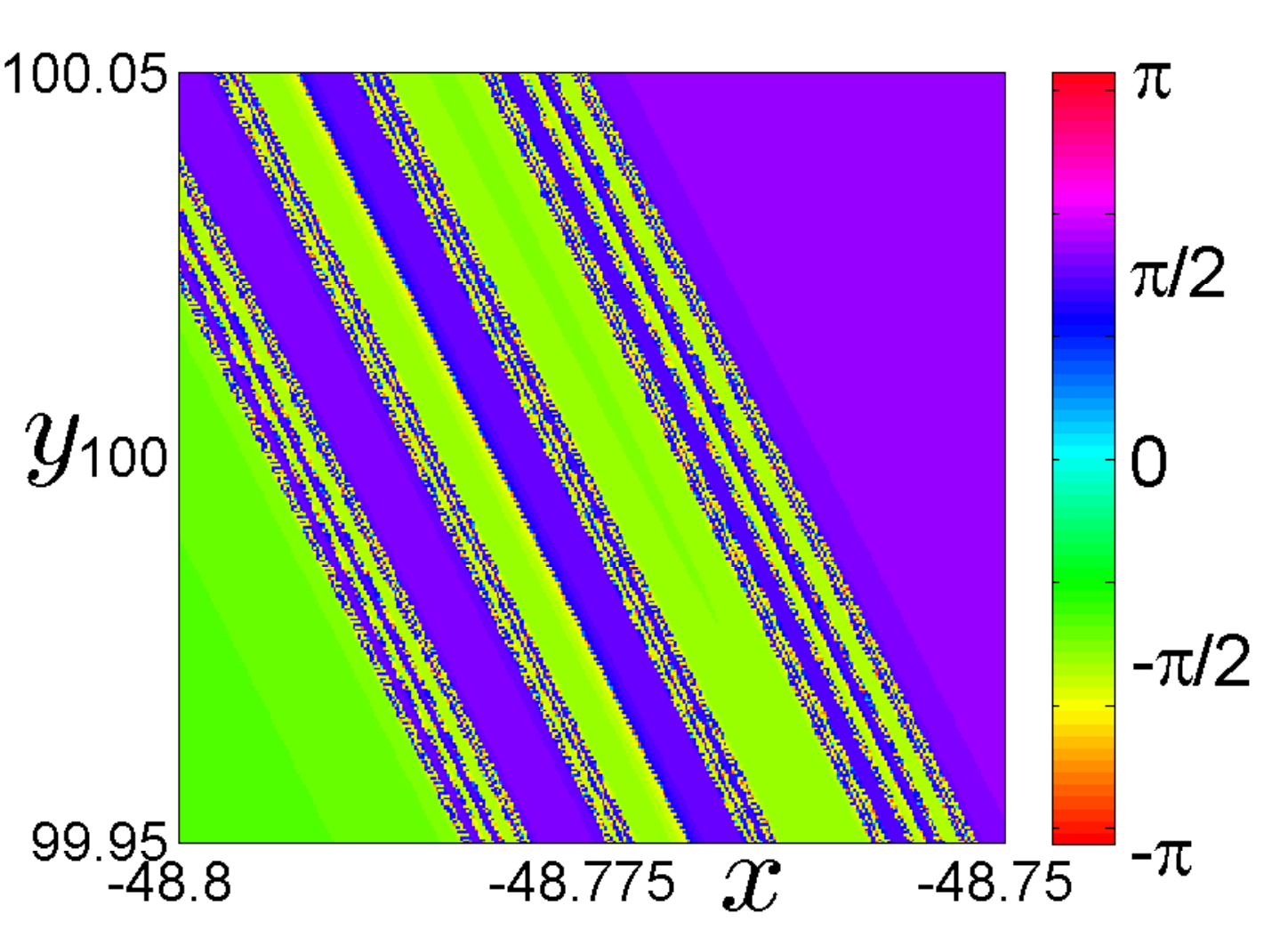}}
\caption{(a)-(c) For the Lorenz model \eqref{Lorenz}, the phase function exhibits complex patterns with boundaries characterized by high phase variation (phaseless set). (The cross-section (c) was chosen so that it contains the two unstable fixed points at $(x,y) \approx \pm(29.17,29.17)$. The black curve is the projection of the limit cycle on the cross-section.) (d) The phaseless set $\mathcal{S}$ has a fractal geometry. (The close-up focuses on a region of the cross-section $z=319$, marked with the black circle in (c).)}
\label{Lorenz_global}
\end{center}
\end{figure}

\red{For the Lorenz system (and in general, for asymptotically periodic systems that exhibit transient chaos), the particular fractal properties of the phaseless set have a significant effect on the phase function, which exhibits unusually complex patterns (Figure \ref{Lorenz_global}(a)-(c)).} While it is well-known that the phase function and the isochrons may be complicated near the phaseless set (e.g. near an unstable fixed point, see \cite{Langfield,Osinga}), the remarkable fact relies here in their fractal properties, which is induced by the fractal geometry of the phaseless set itself (see the close-up in Figure \ref{Lorenz_global}(d)). Note that, for the sake of clarity, we show only the phase function in Figure \ref{Lorenz_global}, the fractal isochrons being mainly concentrated near the phaseless set $\mathcal{S}$.

\begin{figure}[h]
\begin{center}
\includegraphics[width=9cm]{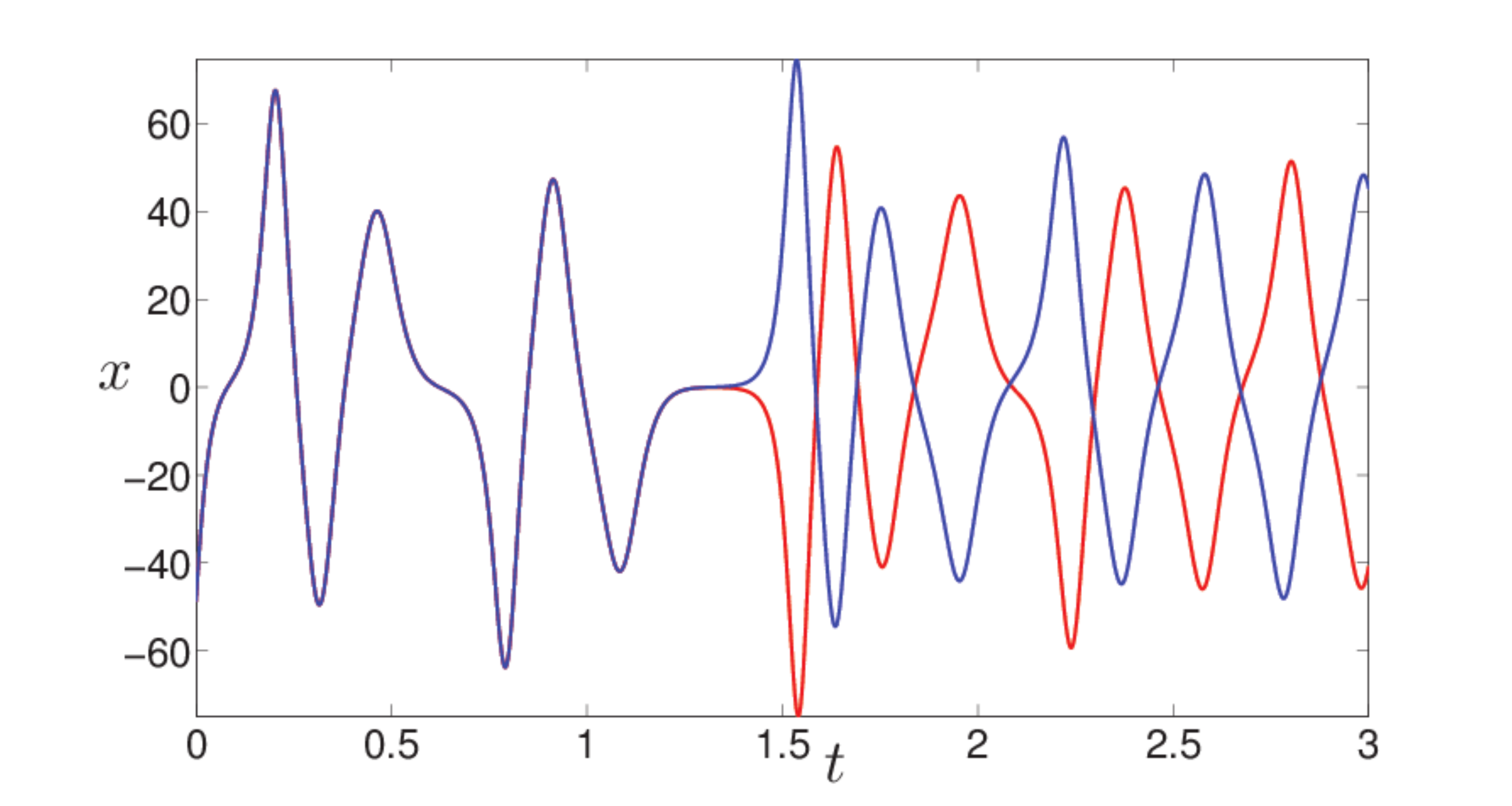}
\caption{As an illustration of the high phase sensitivity observed near the phaseless set, in the Lorenz system \eqref{Lorenz}, two trajectories with close initial conditions diverge (after a long time period) and exhibit two different asymptotic behaviors on the limit cycle. The initial conditions are chosen in the close-up of Figure \ref{Lorenz_global}(d) ($x(0)=-48.7810$ (red trajectory) and $-48.7810-2\times 10^{-4}$ (blue trajectory), $y(0)=100$, $z(0)=319$).}
\label{Lorenz_split}
\end{center}
\end{figure}

In presence of transient chaos, the neighborhood of $\mathcal{S}$ is characterized by a (extremely) high phase sensitivity. Figure \ref{Lorenz_split} shows that two trajectories starting from this region, with very close initial conditions, can have different behaviors reflected in their asymptotic phase. During a time of the order of several limit cycle periods, the two trajectories remain close to the fractal phaseless set $\mathcal{S}$, then \guillemets{escape} $\mathcal{S}$ and diverge near the $z$ axis, subsequently reaching different regions of the limit cycle. Note that it is not sensitivity to initial conditions in classical sense, where exponential \red{divergence} is forever (i.e. positive Lyapunov exponent). Instead, it is the popular notion of phase sensitivity where small changes in initial conditions can separate the trajectories so that they are characterized by different asymptotic behaviors on the limit cycle (i.e. different phases).

\begin{figure}[h]
\begin{center}
\includegraphics[height=4.8cm]{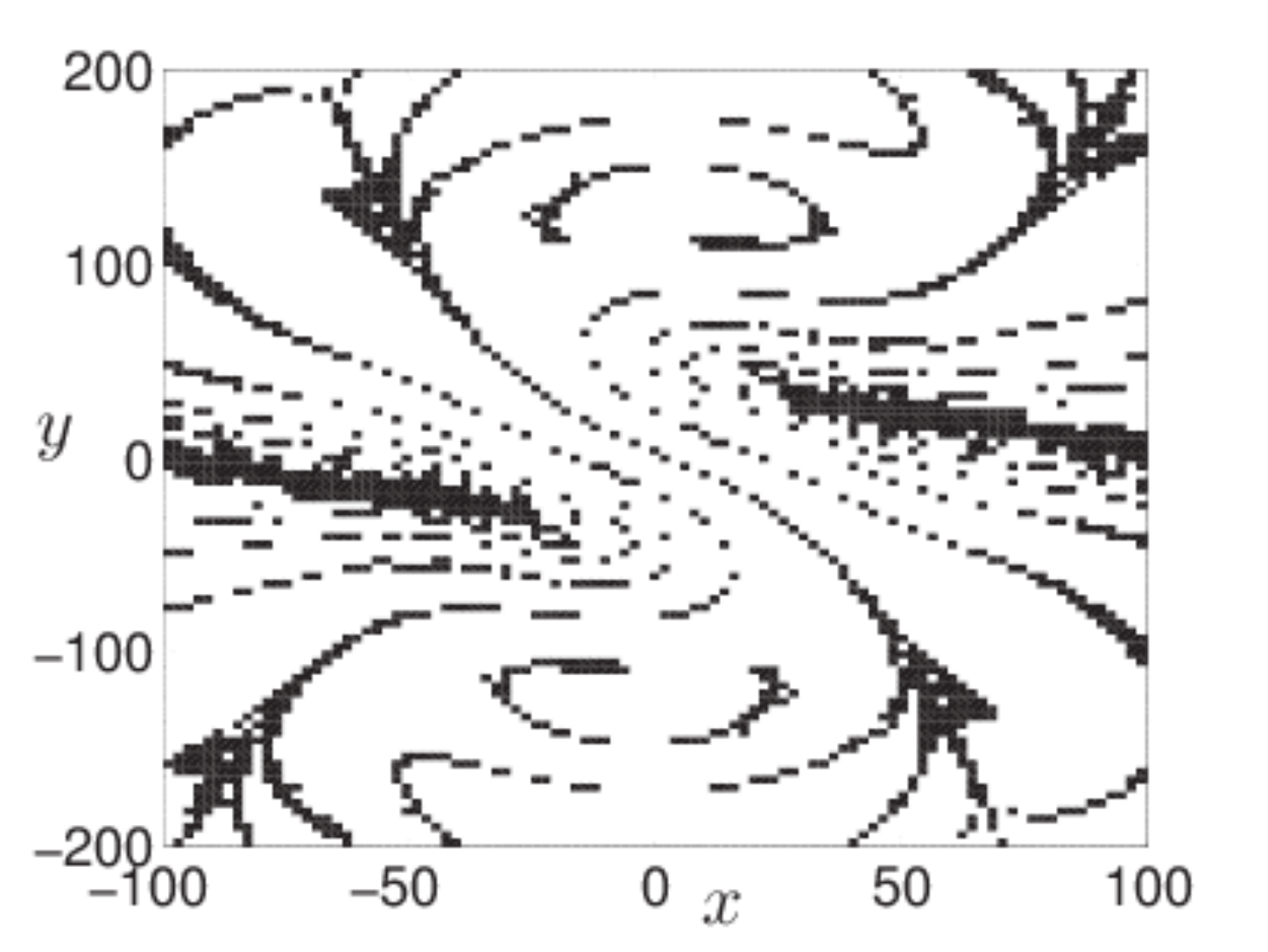}
\caption{Regions of high finite-time Lyapunov exponent (in black) are characterized by a high sensitivity to initial conditions (i.e. high phase sensitivity). For the Lorenz model \eqref{Lorenz} (in the cross-section $z=319$), one verifies that these regions lie in the neighborhood of the fractal phaseless set (see Figure \ref{Lorenz_global}(c)). [The largest finite-time Lyapunov exponent is computed with a finite horizon $T=10$. The black region corresponds to a value higher than $0.25$.]}
\label{Lyap_expo}
\end{center}
\end{figure}

As shown in Figure \ref{Lyap_expo}, this phenomenon can be captured through the computation of the largest finite-time (or local) Lyapunov exponent \cite{local_lyap_expo2,local_lyap_expo1}. For given initial condition $\ve{x}$ and time horizon $T$, the finite-time Lyapunov exponents are given by the logarithm of the eigenvalues of the matrix
\begin{equation*}
\Lambda = \left(M^\textrm{T}(T) M(T) \right)^{1/(2T)}
\end{equation*}
\red{where $M^\mathrm{T}$ denote the transpose of $M$.} For continuous-time systems, $M(\cdot)$ is the fundamental matrix solution of
\begin{equation*}
\frac{d M}{dt}=J(\varphi(t,\ve{x})) M 
\end{equation*}
with $M(0)=I$ and $J$ is the Jacobian matrix of the vector field $\ve{F}$. For discrete-time maps, we have
\begin{equation*}
M(T)=\prod_{t=0}^{T-1} J(\varphi(t,\ve{x}))\,.
\end{equation*}

Regions of high finite-time Lyapunov exponent (black regions) are associated with a high sensitivity to initial conditions. By comparing with Figure \ref{Lorenz_global}(c), we verify that these regions lie close to the fractal phaseless set.


\subsection{Discrete-time model}

\begin{figure}[h!]
\begin{center}
\subfigure[]{\includegraphics[width=7.9cm]{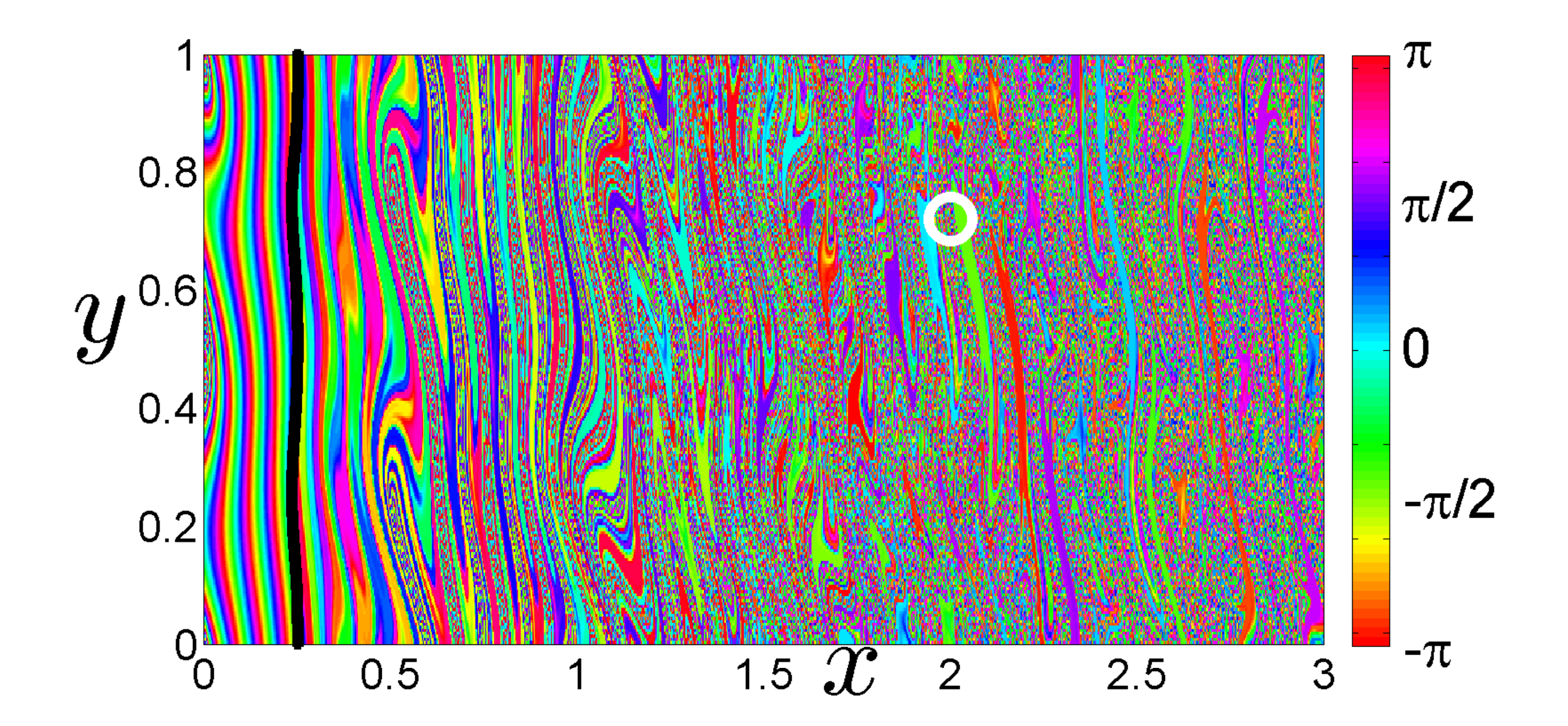}}
\subfigure[]{\includegraphics[width=7.7cm]{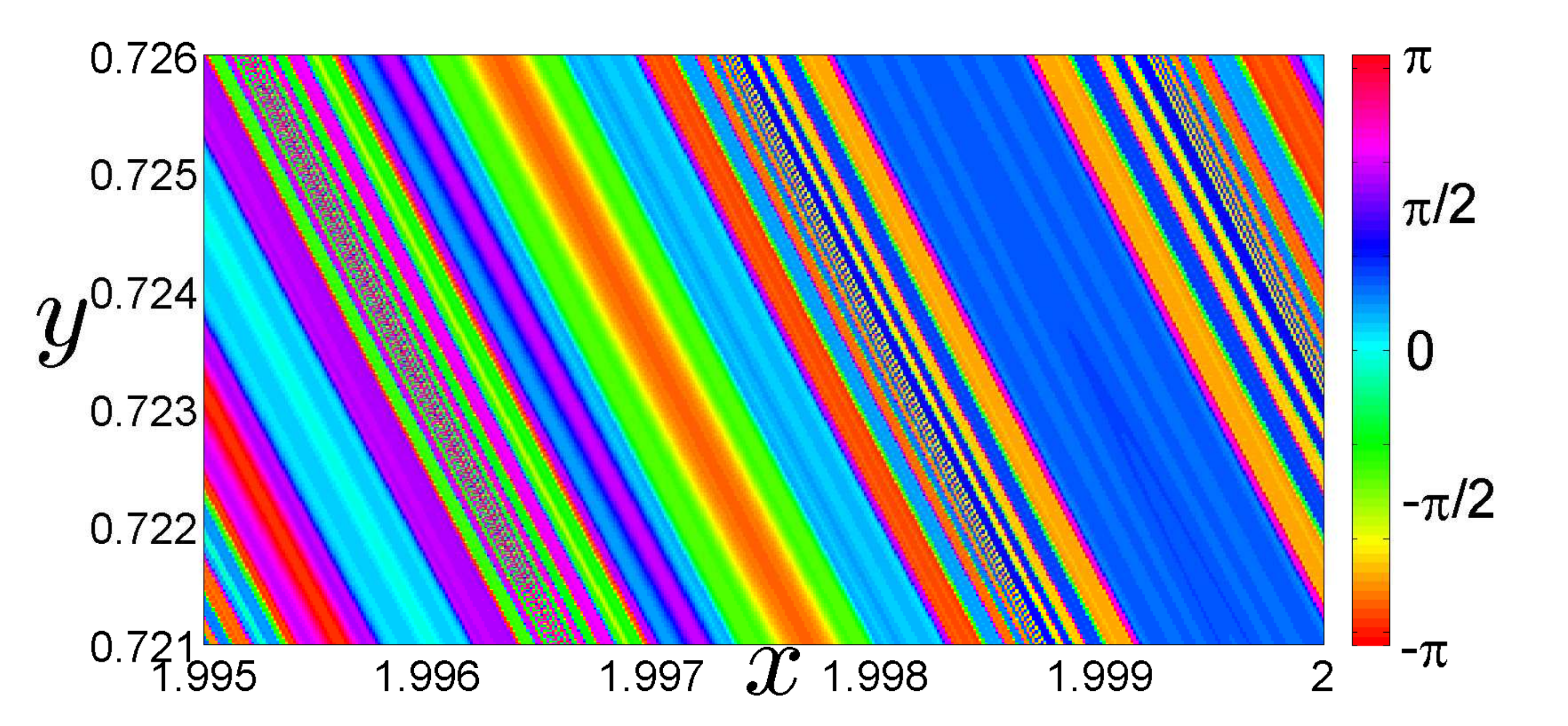}}
\caption{(a) The asymptotic phase associated with the discrete map \eqref{discrete_map} is characterized by complex patterns. The black curve is the \red{attractor $\Gamma$.} (b) A close-up in a region of the state space (marked with the white circle in (a)) shows the self-similarity property of the asymptotic phase.}
\label{discrete_phase}
\end{center}
\end{figure}

\begin{figure}[h]
\begin{center}
\includegraphics[width=10cm]{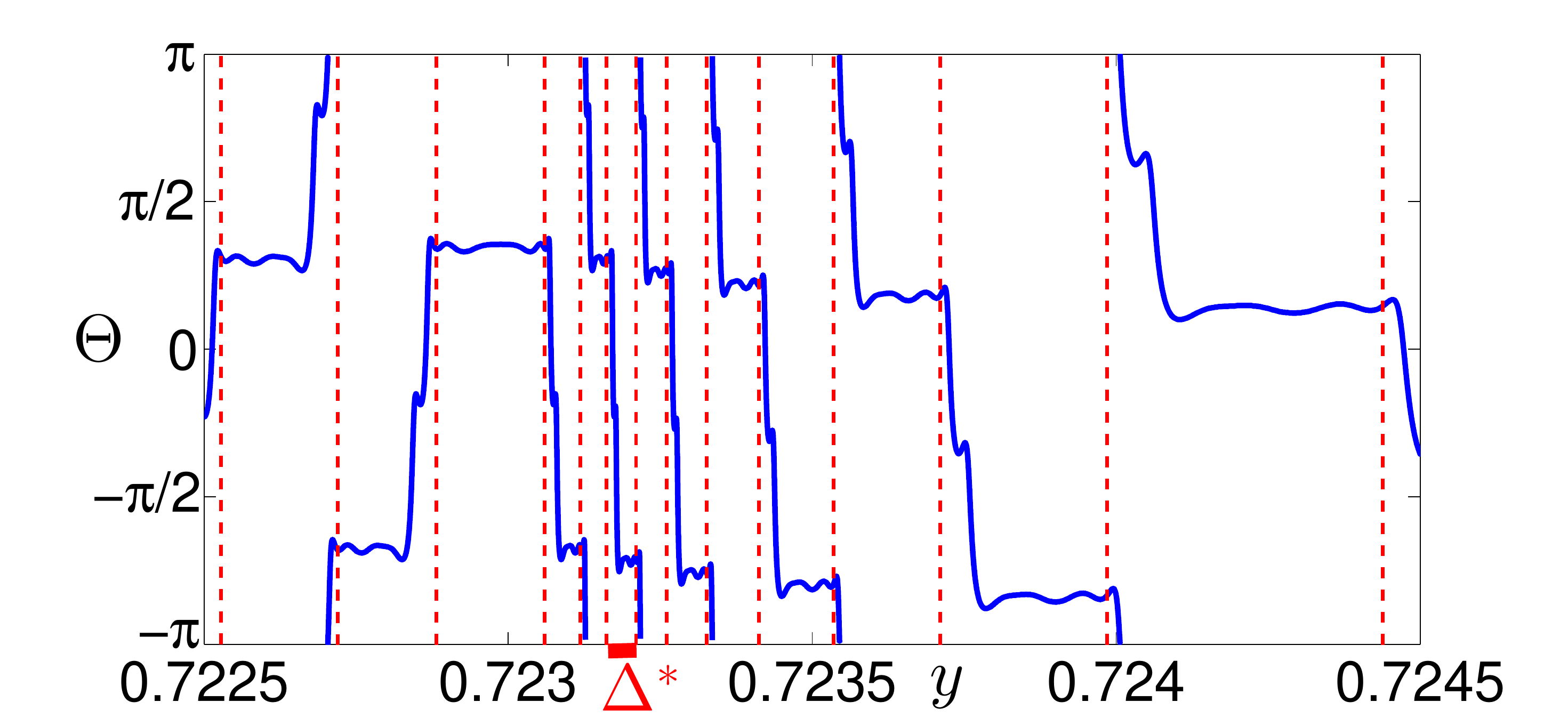}
\caption{For the map \eqref{discrete_map}, the phase function along the line $x=2$ shows that self-similar patterns (separated by the red dashed lines) are not observed at every scale. The size $\Delta^*$ of the smallest pattern is of the order of $3 \times 10^{-5}$. (The resolution of the curve is $6.67 \times 10^{-8}$.)}
\label{discrete_phase_line}
\end{center}
\end{figure}

\begin{figure}[h]
\begin{center}
\includegraphics[height=4.8cm]{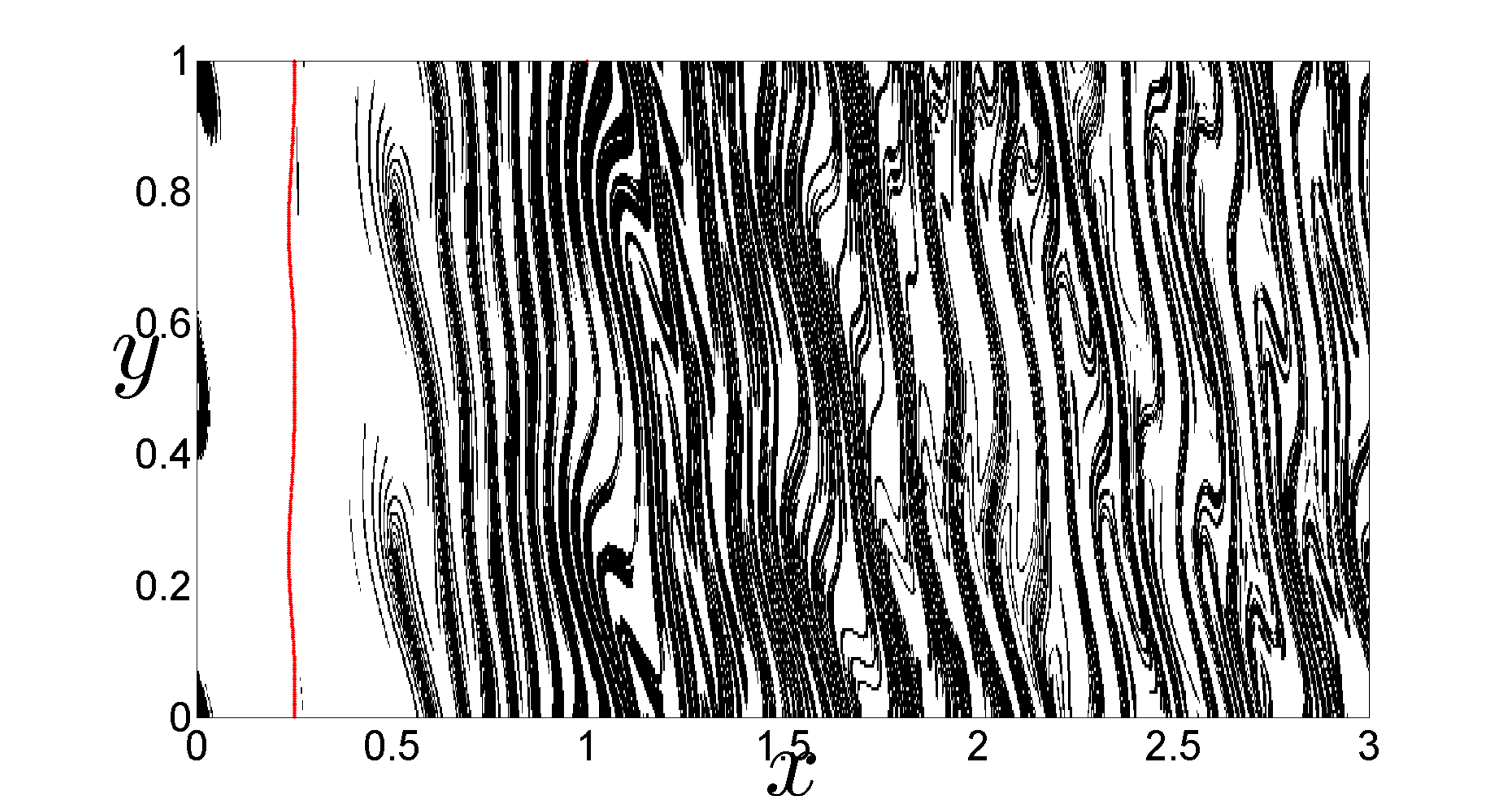}
\caption{For the map \eqref{discrete_map}, regions of high finite-time Lyapunov exponent (in black) are characterized by a high sensitivity to initial conditions (i.e. high phase sensitivity). These regions lie in the neighborhood of the almost phaseless set. The red curve is the invariant dense orbit. [The largest finite-time Lyapunov exponent is computed over $T=35$ iterations of the map. The black region corresponds to a value higher than $0.1$.]}
\label{Lyap_expo2}
\end{center}
\end{figure}

We consider the following map (taken from \cite{MezicBana}):
\begin{equation}
\label{discrete_map}
\begin{array}{rcl}
x(t+1) & = & (1-\gamma) x(t) + a \sin^2(2\pi y(t)) \\
y(t+1) & = & x(t)+y(t)+a \sin(2\pi y(t)) \bmod 1
\end{array}
\end{equation}
with $a=0.03$ and $\gamma=0.06123456756432$. The map has \red{an attractor $\Gamma$ (a topological circle on which the dynamics has a rotation number $\nu_0 \approx 0.24482525$)} near $x=0.25$. The corresponding asymptotic phase function is characterized by a complex geometry (see Figure \ref{discrete_phase}(a) or Figure 12(a) in \cite{MezicBana}) with self-similar patterns (Figure \ref{discrete_phase}(b)). In contrast to the continuous-time Lorenz system, these patterns are not observed at very small scales, but have a minimal size $\Delta^*$ that depends on the distance to $\Gamma$ (Figure \ref{discrete_phase_line}). (Note that the patterns have an arbitrarily small size as they are observed far from $\Gamma$, i.e. for $x\rightarrow \infty$.) The region of self-similar patterns is therefore an \emph{almost phaseless set}\footnote{\red{The term \guillemets{almost phaseless set} has been coined in \cite{Osinga}}} $\tilde{\mathcal{S}}$ rather than a phaseless set, i.e. a region of very large---but not arbitrarily large---phase variation. It is characterized by a high sensitivity to initial conditions, as shown by the values of the largest finite-time Lyapunov exponent (Figure \ref{Lyap_expo2}). \red{It is important to note that there is no precise definition of the almost phaseless set, which depends on the threshold chosen to assess whether the phase variation is large or not. Moreover, the almost phaseless set does not correspond to a specific invariant set known a priori.} In contrast to \red{the phaseless set in the continuous-time situation}, it cannot be interpreted in terms of a fractal chaotic saddle, since all trajectories---including those with the initial conditions in $\tilde{\mathcal{S}}$---converge to \red{the attractor $\Gamma$.}

\newpage
\section{Phase sensitivity and fractal dimension}
\label{sec_sensitivity_fractal_dim}

Phaseless sets are always characterized by a high sensitivity of the asymptotic phase. But when they have a fractal geometry, they occupy an important portion of the state space, so that the overall phase sensitivity of the system is accentuated. The largest  Lyapunov exponent cannot capture this overall phase sensitivity, since it is always equal to zero (the finite-time Lyapunov exponent tends to zero as the time horizon increases). Also, although the computation of finite-time Lyapunov exponents can be used to highlight sensitive regions of the state space (see Figures \ref{Lyap_expo} and \ref{Lyap_expo2}), it only provides a local measure of the system sensitivity which depends on the chosen time horizon and initial condition. Therefore, a new notion is required to measure the overall phase sensitivity of the system. To this end, we define a \emph{phase sensitivity coefficient} and we show that this coefficient is closely related to the fractal dimension of the isochrons. This implies that the coefficient is \red{an invariant of the system} and can be used to compare the phase sensitivity of different systems.

\subsection{Phase sensitivity}

Consider the geodesic distance $d:\mathbb{S}^1 \times \mathbb{S}^1 \rightarrow [0,\pi)$ on the circle
\begin{equation}
\label{geo_dist}
d(\theta, \theta')=\min_{k\in \mathbb{Z}} |\theta-\theta'+k 2 \pi|\,.
\end{equation}
We define the \emph{phase sensitivity function} $f:\mathcal{B} \times \mathbb{R}^+ \rightarrow [0,\pi)$ by
\begin{equation}
\label{phase_sens_fct}
f(\ve{x},\epsilon)=\max_{\ve{x}'\in B(\ve{x},\epsilon) \cap \mathcal{B}} d(\Theta(\ve{x}),\Theta(\ve{x}'))\,,
\end{equation}
where $\mathcal{B}$ is the basin of attraction of the limit cycle and $B(\ve{x},\epsilon)$ is a ball with center $\ve{x}$ and radius $\epsilon>0$. If there is an uncertainty $\epsilon$ (e.g. induced by external perturbations or noise) on the initial condition $\ve{x}$ of a trajectory, then the asymptotic behavior of the trajectory will be associated with an uncertainty $f(\ve{x},\epsilon)$ on the phase. Note that $f(\ve{x},\epsilon)$ is the worst-case uncertainty.

\begin{figure}[h]
\begin{center}
\includegraphics[width=8cm]{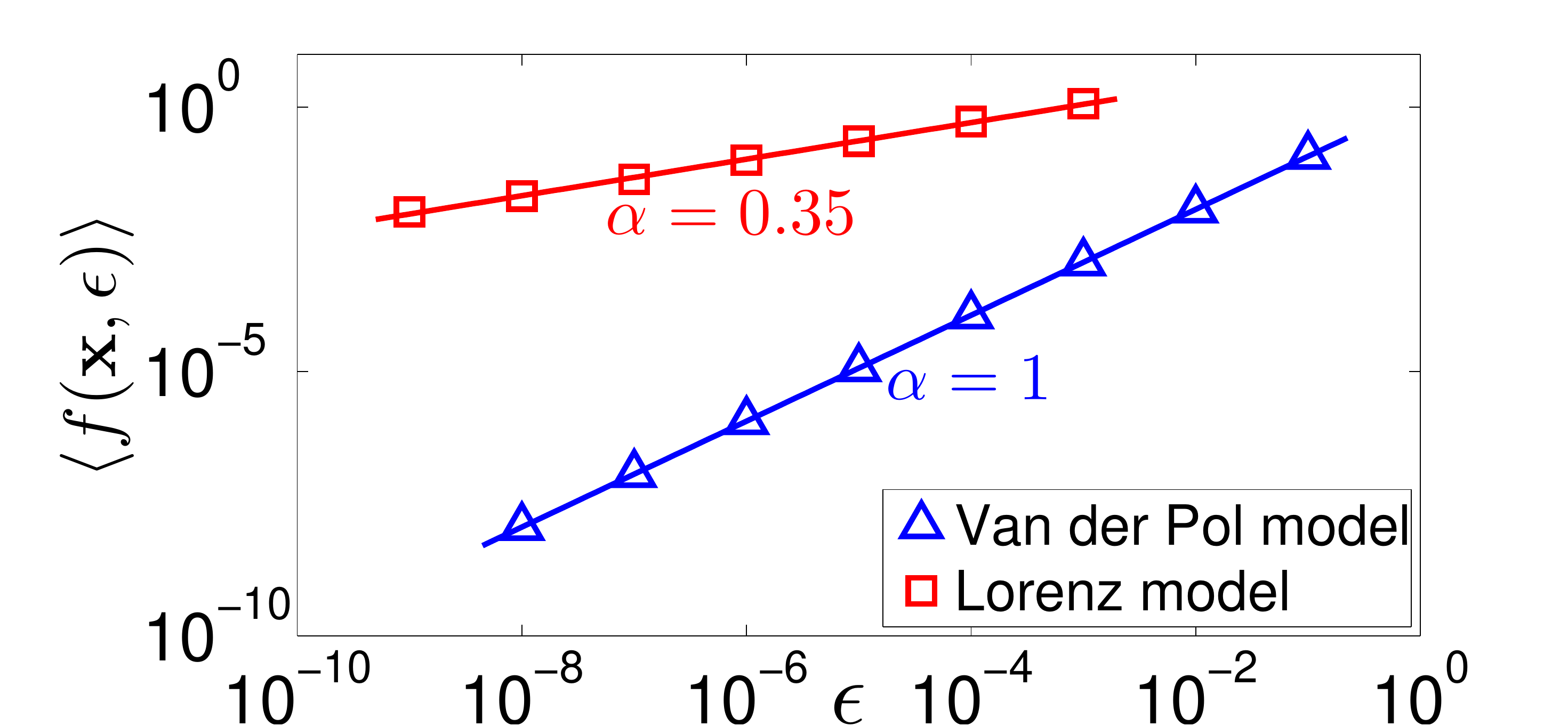}
\caption{When the phase function exhibits fractal patterns, the phase sensitivity coefficient $\beta=1-\alpha$ is greater than $0$. We obtain $\beta=0$ for the Van der Pol model (see Figure \ref{vdp_iso}), and $\beta=0.65$ for the Lorenz model \eqref{Lorenz}. Simulation details are given in Appendix \ref{app_numerical}.}
\label{phase_uncertainty}
\end{center}
\end{figure}

For a given compact set $\mathcal{A} \subset \mathbb{R}^n$, the average phase sensitivity function is computed as
\begin{equation*}
\langle f(\ve{x},\epsilon) \rangle_{\mathcal{A} \cap \mathcal{B}}= \frac{1}{\mu[\mathcal{A} \cap \mathcal{B}]} \int_{\mathcal{A} \cap \mathcal{B}} f(\ve{x},\epsilon) \, d\ve{x}\,,
\end{equation*}
where $\mu$ is the Lebesgue measure on $\mathcal{A}$. \red{If $\mathcal{A} \cap \mathcal{S} \neq \emptyset$ (unless $\mathcal{S} = \emptyset$), we define}
\begin{equation}
\label{phase_uncert_coeff}
\alpha=\lim_{\epsilon \rightarrow 0} \frac{\ln \langle f(\ve{x},\epsilon) \rangle_{\mathcal{A} \cap \mathcal{B}}}{\ln \epsilon}\qquad \beta=1-\lim_{\epsilon \rightarrow 0} \frac{\ln \langle f(\ve{x},\epsilon) \rangle_{\mathcal{A} \cap \mathcal{B}}}{\ln \epsilon}\,,
\end{equation}
or equivalently, $\langle f(\ve{x},\epsilon) \rangle_{\mathcal{A} \cap \mathcal{B}} \sim \epsilon^\alpha=\epsilon^{1-\beta}$ for $\epsilon \ll 1$. We refer to $\beta=1-\alpha$ as the \emph{phase sensitivity coefficient}. If $\beta=0$, reducing the uncertainty $\epsilon$ on the initial condition by a certain amount (e.g. reducing the noise intensity) reduces the average phase uncertainty by the same amount (at least when $\epsilon$ is small). This is the usual situation observed with globally asymptotically stable periodic systems (e.g. Van der Pol model, see Figure \ref{vdp_iso}). But if the phase function has fractal properties, we observe that $\beta>0$. In this case, a reduction of the uncertainty $\epsilon$ produces only a slight reduction of the average phase uncertainty. For instance, in the Lorenz model considered in Section \ref{sec_continuous_model}, a logarithmic plot of the average phase sensitivity function with respect to $\epsilon$ shows that the phase sensitivity coefficient is equal to $\beta \approx 0.65$ (see Figure \ref{phase_uncertainty}). In this case, a reduction of the uncertainty $\epsilon$ by a factor $2$ only reduces the average phase uncertainty by a factor $2^{1-0.65} \approx 1.27$. The phase sensitivity coefficient is therefore directly related to the overall phase sensitivity of the system.

\subsection{Fractal dimension}
\label{sub_sec_fractal_dim}

We now show that the phase sensitivity coefficient \eqref{phase_uncert_coeff} is an invariant of the system that is closely related to the fractal co-dimension of the phaseless set and the isochrons. For a given $\epsilon>0$, consider the set
\begin{equation*}
\mathcal{M}_\mathcal{S}(\epsilon) = \{\ve{x} \in \mathcal{A} \cap \mathcal{B} | B(\ve{x},\epsilon) \cap \mathcal{S} \neq \emptyset \}\,,
\end{equation*}
i.e. the sets of points lying within a distance $\epsilon$ of the phaseless set $\mathcal{S}$. Since the values of the phase function span $\mathbb{S}^1$ on any neighborhood of $\mathcal{S}$, we have
\begin{equation}
\label{integral1}
f(\ve{x},\epsilon)=\pi \qquad \forall \ve{x} \in \mathcal{M}_\mathcal{S}(\epsilon) \,.
\end{equation}
If the limit cycle $\Gamma$ is normally hyperbolic, the isochrons are as smooth as the vector field, and so is the phase function $\Theta$ in $\mathcal{B}$, since $\Theta$ smoothly increases from one isochron to another. It follows that, for $\epsilon$ small enough and $\delta \epsilon>0$, we have
\begin{equation}
\label{integral2}
f(\ve{x},\epsilon)=\max_{\|\ve{e}\|=\epsilon} |\nabla \Theta(\ve{x})\cdot \ve{e}| + \mathcal{O}(\epsilon^2)= \epsilon \|\nabla \Theta(\ve{x})\| + \mathcal{O}(\epsilon^2) \qquad \forall \ve{x} \in \mathcal{A} \cap \mathcal{B} \setminus \mathcal{M_S}(\epsilon+ \delta \epsilon)\,,
\end{equation}
where $\nabla \Theta$ is the gradient of the phase function \red{and $\|\cdot\|$ is the Euclidean norm.} Then \eqref{integral1} and \eqref{integral2} imply that
\begin{equation*}
\begin{split}
\langle f(\ve{x},\epsilon) \rangle_{\mathcal{A} \cap \mathcal{B}} & = \frac{1}{\mu[\mathcal{A}\cap\mathcal{B}]} \left(\int_{\mathcal{M_S}(\epsilon)} f(\ve{x},\epsilon) \, d\ve{x} + \int_{\mathcal{A} \cap \mathcal{B} \setminus \mathcal{M_S}(\epsilon+\delta \epsilon)} f(\ve{x},\epsilon) \, d\ve{x} + \int_{\mathcal{M_S}(\epsilon+ \delta \epsilon) \setminus \mathcal{M_S}(\epsilon)} f(\ve{x},\epsilon) \, d\ve{x} \right) \\
& = \frac{\mu[\mathcal{M}_\mathcal{S}(\epsilon)]}{\mu[\mathcal{A}\cap\mathcal{B}]} \pi + \left(1-\frac{\mu[\mathcal{M}_\mathcal{S}(\epsilon+\delta \epsilon)]}{\mu[\mathcal{A}\cap\mathcal{B}]}\right) \, \epsilon \langle \|\nabla \Theta\| \rangle_{\mathcal{A} \cap \mathcal{B} \setminus \mathcal{M}_\mathcal{S}(\epsilon+\delta \epsilon)} \\
& \qquad \qquad + \frac{1}{\mu[\mathcal{A}\cap\mathcal{B}]} \int_{\mathcal{M_S}(\epsilon+ \delta \epsilon) \setminus \mathcal{M_S}(\epsilon)} f(\ve{x},\epsilon) \, d\ve{x} + \mathcal{O}(\epsilon^2)\,.
\end{split}
\end{equation*}
Since by definition $0 \leq f(\ve{x},\epsilon) \leq \pi$, one has
\begin{equation*}
0 \leq \int_{\mathcal{M_S}(\epsilon+ \delta \epsilon) \setminus \mathcal{M_S}(\epsilon)} f(\ve{x},\epsilon) \, d\ve{x} \leq \mu[\mathcal{M_S}(\epsilon+ \delta \epsilon) \setminus \mathcal{M_S}(\epsilon)]  \pi
\end{equation*}
and it follows that
\begin{equation*}
\begin{split}
& \frac{\mu[\mathcal{M}_\mathcal{S}(\epsilon)]}{\mu[\mathcal{A}\cap\mathcal{B}]} \pi + \left(1-\frac{\mu[\mathcal{M}_\mathcal{S}(\epsilon+\delta \epsilon)]}{\mu[\mathcal{A}\cap\mathcal{B}]}\right) \, \epsilon \langle \|\nabla \Theta\| \rangle_{\mathcal{A} \cap \mathcal{B} \setminus \mathcal{M}_\mathcal{S}(\epsilon+\delta \epsilon)} + \mathcal{O}(\epsilon^2) \\
&  \leq \langle f(\ve{x},\epsilon) \rangle_{\mathcal{A} \cap \mathcal{B}} \leq \frac{\mu[\mathcal{M}_\mathcal{S}(\epsilon+\delta \epsilon)]}{\mu[\mathcal{A}\cap\mathcal{B}]} \pi + \left(1-\frac{\mu[\mathcal{M}_\mathcal{S}(\epsilon+\delta \epsilon)]}{\mu[\mathcal{A}\cap\mathcal{B}]}\right) \, \epsilon \langle \|\nabla \Theta\| \rangle_{\mathcal{A} \cap \mathcal{B} \setminus \mathcal{M}_\mathcal{S}(\epsilon+\delta \epsilon)} + \mathcal{O}(\epsilon^2) \,.
\end{split}
\end{equation*}
Considering $\delta \epsilon = \epsilon$ and taking the limit $\epsilon \rightarrow 0$ yield

\begin{equation}
\label{ratio_uncert_coeff}
\begin{split}
& \lim_{\epsilon \rightarrow 0} \frac{\ln \left( \mu[\mathcal{M}_\mathcal{S}(\epsilon)] + (1-\mu[\mathcal{M}_\mathcal{S}(2\epsilon)]/\mu[\mathcal{A}\cap\mathcal{B}]) \, \epsilon \langle \|\nabla \Theta\| \rangle_{\mathcal{A} \cap \mathcal{B}\setminus \mathcal{M}_\mathcal{S}(2 \epsilon)}\right)}{\ln \epsilon} \\
& \geq \alpha \geq \lim_{\epsilon \rightarrow 0} \frac{\ln \left( \mu[\mathcal{M}_\mathcal{S}(2\epsilon)] + (1-\mu[\mathcal{M}_\mathcal{S}(2\epsilon)]/\mu[\mathcal{A}\cap\mathcal{B}]) \, \epsilon \langle \|\nabla \Theta\| \rangle_{\mathcal{A} \cap \mathcal{B}\setminus \mathcal{M}_\mathcal{S}(2 \epsilon)}\right)}{\ln \epsilon} \,,
\end{split}
\end{equation}
where we used $\lim_{\epsilon \rightarrow 0} \ln (\pi/\mu[\mathcal{A}\cap\mathcal{B}])/\ln \epsilon=0$ and the fact that $\ln \epsilon <0$ for $\epsilon \ll 1$.
If $\mathcal{M}_\mathcal{S}(\epsilon) = \emptyset$ (i.e. $\mathcal{S}=\emptyset$), we have simply
\begin{equation*}
\alpha=\lim_{\epsilon \rightarrow 0} \frac{\ln \left(\epsilon \langle \|\nabla \Theta\| \rangle_{\mathcal{A} \cap \mathcal{B}} \right)}{\ln \epsilon} = 1 \qquad \beta=0\,.
\end{equation*}
Otherwise, according to \cite{McDonald}, $\mu[\mathcal{M}_\mathcal{S}(\epsilon)]$ scales as
\begin{equation}
\label{scaling_dim}
\mu[\mathcal{M}_\mathcal{S}(\epsilon)] \sim \epsilon^{N-D} \,,
\end{equation}
with $N$ the dimension of the state space and $D$ the capacity dimension of the fractal set $\mathcal{S}$, \red{i.e. (see \cite{Farmer_fractal_dim})
\begin{equation}
\label{def_capacity_dimension}
D=\lim_{d \rightarrow 0} \frac{\ln \mathcal{N}(d)}{\ln(1/d)}
\end{equation}
where $\mathcal{N}(d)$ is the number of boxes of size $d$ required to cover $\mathcal{S}$.} In addition, if $\mathcal{S}$ is a normally hyperbolic repeller of co-dimension (at most) one, the gradient $\|\nabla \Theta\|$ scales as $1/\epsilon$ when it is evaluated at a small distance $\epsilon$ of $\mathcal{S}$ (see Appendix \ref{app_scaling_gradient}). \red{It follows that its integral scales as $\ln \epsilon$ when the boundary of the domain of integration is at a distance $\epsilon$ of $\mathcal{S}$}, so that $\langle \|\nabla \Theta\| \rangle_{\mathcal{A} \cap \mathcal{B}\setminus \mathcal{M}_\mathcal{S}(2 \epsilon)} \sim \ln \epsilon$ for $\epsilon \ll 1$ \footnote{If $\mathcal{S}$ is of co-dimension greater than one, intuitive arguments suggest that $\langle \|\nabla \Theta\| \rangle_{\mathcal{A} \cap \mathcal{B}\setminus \mathcal{M}_\mathcal{S}(2 \epsilon)} = \mathcal{O}(1)$. However, a rigorous result is beyond the scope of this paper.}. Finally, \eqref{ratio_uncert_coeff} leads to
\begin{equation*}
\lim_{\epsilon \rightarrow 0} \frac{\ln \left( \epsilon^{N-D} + (1-(2\epsilon)^{N-D}) \, \epsilon \ln \epsilon \right)}{\ln \epsilon} \geq \alpha \geq \lim_{\epsilon \rightarrow 0} \frac{\ln \left( (2\epsilon)^{N-D} + (1-(2\epsilon)^{N-D}) \, \epsilon \ln \epsilon\right)}{\ln \epsilon}
\end{equation*}
or equivalently
\begin{equation}
\label{rel_fractal_dim}
\alpha=N-D \qquad \beta=1-(N-D)
\end{equation}
(for $\mathcal{S}$ of co-dimension $N-D \leq 1$).
This important relationship implies that the phase sensitivity coefficient is directly related to the fractal dimension of $\mathcal{S}$. More precisely, the coefficient will be strictly greater than $0$ only if $\mathcal{S}$ is fractal. (Note that a phaseless set $\mathcal{S}$ of co-dimension greater than $1$ is associated with a phase sensitivity coefficient equal to $0$, even if it is fractal.) This result also implies that the phase sensitivity coefficient is an intrinsic property of the system. It has a unique value that does not depend on the choice of the set $\mathcal{A}$.

The fractal properties of $\mathcal{S}$ and the phase sensitivity can also be characterized by considering the set
\begin{equation*}
\mathcal{M}_{\delta \theta}(\epsilon)=\{\ve{x} \in \mathcal{A} \cap \mathcal{B} |\exists \ve{x}' \in B(\ve{x},\epsilon) \textrm{ s.t. } d(\Theta(\ve{x}),\Theta(\ve{x}'))>\delta \theta \}
\end{equation*}
for given $\epsilon>0$ and $0<\delta \theta<\pi$. Since it is clear that $\mathcal{M}_{\delta \theta}(\epsilon) \approx \mathcal{M}_\mathcal{S}(\epsilon)$ as $\epsilon\rightarrow 0$, \eqref{scaling_dim} yields the additional relationship
\begin{equation}
\label{sensitivity}
N-D = \lim_{\epsilon \rightarrow 0} \frac{\ln \mu [\mathcal{M}_{\delta \theta}(\epsilon)]}{\ln \epsilon}
\end{equation}
for any value $\delta \theta$, provided that $\mathcal{M}_{\delta \theta}(\epsilon, \delta \theta) \neq \emptyset$. If $\mathcal{S}$ is fractal with $N-D<1$, then the fraction of trajectories for which an uncertainty smaller than $\epsilon$ on the initial condition induces a phase uncertainty greater than $\delta \theta$ on the asymptotic phase will only be slightly reduced by a significant decreasing of $\epsilon$. We remark that this property is independent of the value $\delta \theta$.

\begin{figure}[h]
\begin{center}
\includegraphics[width=8cm]{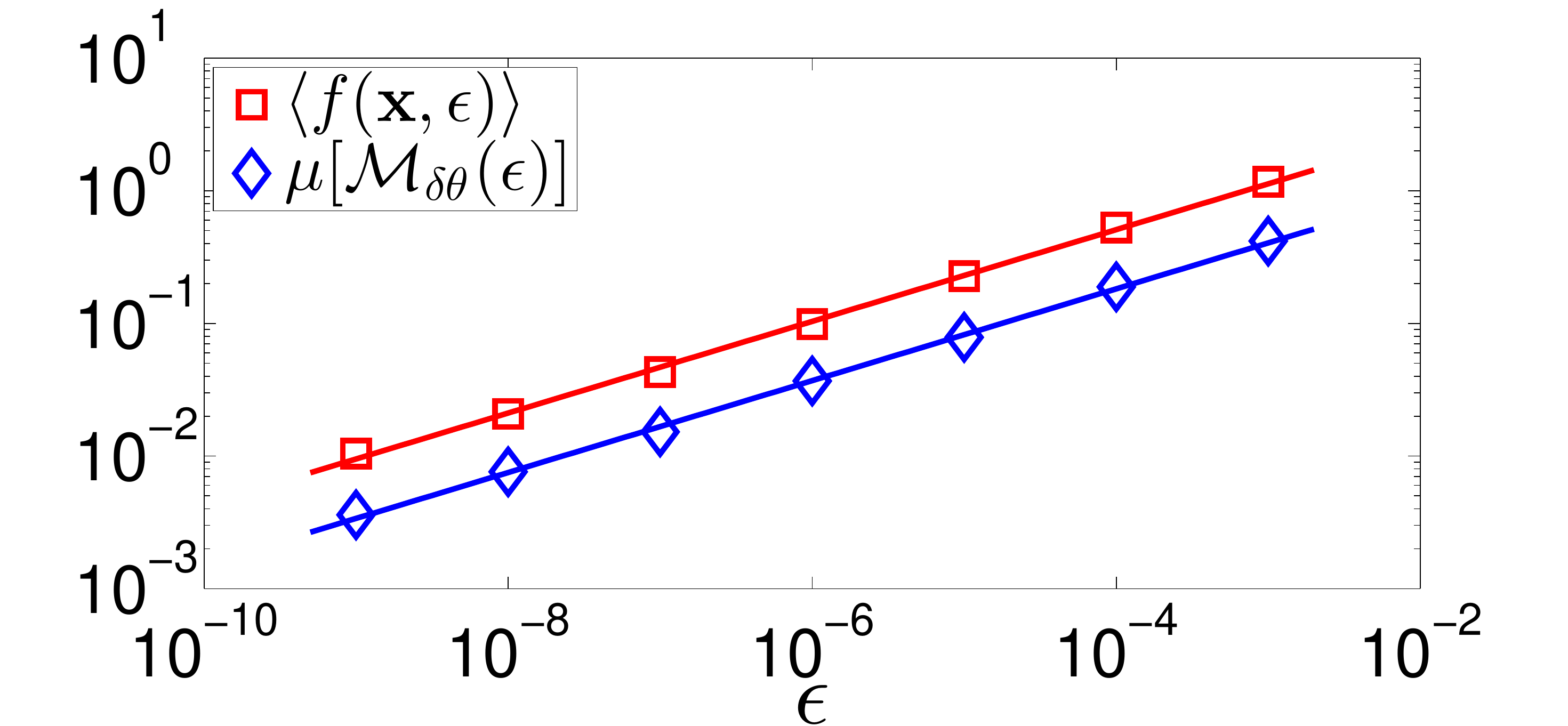}
\caption{The phase sensitivity coefficient and the fractal dimension of $\mathcal{S}$ (and of the isochrons) are computed for the Lorenz model \eqref{Lorenz}, with \eqref{phase_uncert_coeff} (in red) and with \eqref{sensitivity} for the value $\delta \theta =0.5$ (in blue). The two sets of formulae yield equivalent results, both showing that the fractal dimension is equal to $2.65$ ($\alpha=0.34573$ with the first method; $\alpha=0.34623$ with the second method). Simulation details are given in Appendix \ref{app_numerical}.}
\label{fractal_dim_comparison}
\end{center}
\end{figure}

The definitions and equalities \eqref{phase_uncert_coeff}-\eqref{rel_fractal_dim}-\eqref{sensitivity} summarize the relationships between the overall phase sensitivity of the system and the fractal dimension of the isochrons and phaseless set. Figure \ref{fractal_dim_comparison} shows that \eqref{phase_uncert_coeff} and \eqref{rel_fractal_dim}-\eqref{sensitivity} yield similar results.

\paragraph*{The isochrons are fractal.}When the phaseless set $\mathcal{S}$ is fractal, the isochrons themselves are also fractal. Since any neighborhood of $\mathcal{S}$ intersects an isochron $\mathcal{I}_\theta$ \cite{Guckenheimer_iso}, every box of size $d$ used to cover $\mathcal{S}$ intersects $\mathcal{I}_\theta$. Hence, we have
\begin{equation*}
\mathcal{N}_\mathcal{I}(d) \geq \mathcal{N}_\mathcal{S}(d)\,,
\end{equation*}
where $\mathcal{N}_\mathcal{S}(d)$ and $\mathcal{N}_\mathcal{I}(d)$ are the number of boxes of size $d$ required to cover $\mathcal{S}$ and $\mathcal{I}_\theta$, respectively. It follows that, for $d<1$,
\begin{equation*}
\frac{N_\mathcal{I}(d)}{\log(1/d)} \geq \frac{N_\mathcal{S}(d)}{\log(1/d)}
\end{equation*}
and taking the limit $d \rightarrow 0$, we obtain from \eqref{def_capacity_dimension} that the capacity dimension of $\mathcal{I}_\theta$ is greater or equal to the capacity dimension $D$ of $\mathcal{S}$. 

\paragraph*{The phase response curve is fractal.}The (finite) phase response curve \eqref{PRC} is fractal provided that the curve $\Omega=\{\ve{x}^\gamma(\theta)+\ve{e}|\theta \in \mathbb{S}^1\}$ satisfies $\Omega \cap \mathcal{S} \neq \emptyset$. The fractal dimension is obtained as follows. We assume that the isochrons---i.e. the level sets of $\Theta$---have a fractal dimensional equal to $D=N-1+\beta$ (see \eqref{rel_fractal_dim}). Then the hypersurface $\Lambda=\{(\ve{x},\Theta(\ve{x}))|\ve{x}\in \mathcal{B}\} \subset \mathcal{B}\times \mathbb{S}^1$ is of dimension $N+\beta$ and it follows that the curve $\{(\ve{x},\Theta(\ve{x}))|\ve{x} \in \Omega \} \subset \Lambda$ is of dimension $N+\beta-(N-1)=1+\beta$. This implies that the dimension of the phase response curve \eqref{PRC} is also equal to $1+\beta$. An example of fractal phase response curve will be given in Section \ref{sec_bursting} (Figure \ref{fractal_PRC_elliptic}).

\begin{figure}[h]
\begin{center}
\includegraphics[width=7cm]{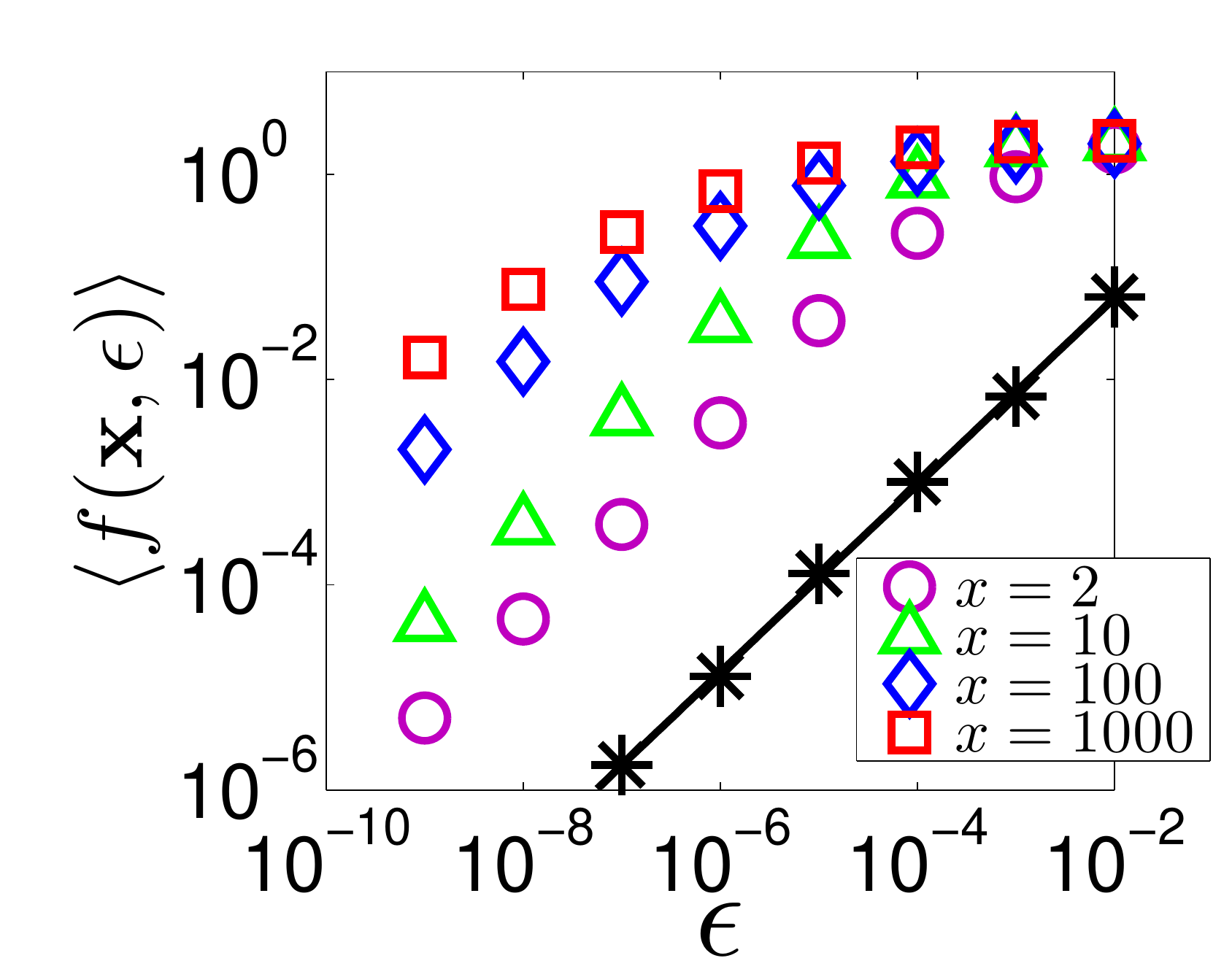}
\caption{For the discrete-time map \eqref{discrete_map}, the phase sensitivity coefficient $\beta=1-\alpha$ is computed on one-dimensional intervals $y\in[0,1]$ with $x\in\{2,10,100,1000\}$ constant. It is always equal to zero, as indicated by the slope of $\ln \langle f \rangle$ which is close to $1$ for small values $\epsilon$. However, for large values $\epsilon$ (i.e. $\epsilon \gg \Delta^*$), the slope is less than one, so that the infinitesimal phase sensitivity coefficient \eqref{inf_phase_sens_coeff} is greater than zero. The phenomenon is observed on a broad interval when the set $\mathcal{A}$ is chosen far from the \red{attractor} (i.e. $x\gg 1$), where $\Delta^*$ is very small. In contrast, for the model of Figure \ref{vdp_iso}(b), there is no fractal (almost) phaseless set and we observe that $\langle f \rangle \sim \epsilon^1$ for all values $\epsilon$ (black stars). Simulation details are given in Appendix \ref{app_numerical}.}
\label{fractal_uncertainty_discrete}
\end{center}
\end{figure}

\paragraph*{Almost phaseless set.} The discrete map \eqref{discrete_map} has an almost phaseless set $\tilde{\mathcal{S}}$ characterized by self-similar patterns scaling down to $\Delta^*$, at best. This implies that the phase sensitivity coefficient, as defined in \eqref{phase_uncert_coeff}, is equal to zero, so that the dimension $D$ of $\tilde{\mathcal{S}}$ is one. But for large values $\epsilon$ (i.e. $\epsilon \gg \Delta^*$), the phase sensitivity function behaves as if the fractal dimension were more than one: the average phase sensitivity function remains close to one and slowly decreases as $\epsilon$ decreases (see Figure \ref{fractal_uncertainty_discrete}). In other words, the \guillemets{infinitesimal} phase sensitivity coefficient
\begin{equation}
\label{inf_phase_sens_coeff}
1-\frac{d }{d (\ln \epsilon)} \ln \langle f(\ve{x},\epsilon) \rangle_{\mathcal{A} \cap \mathcal{B}} \red{ \,= 1-\epsilon \frac{d }{d \epsilon} \ln \langle f(\ve{x},\epsilon) \rangle_{\mathcal{A} \cap \mathcal{B}}}
\end{equation}
\red{where $d/dX$ denotes the derivative with respect to $X$}, is large for large $\epsilon$. In contrast, for maps that do not have an almost phaseless set with these self-similarity properties, the average phase sensitivity function scales as $\epsilon^1$ for all values of $\epsilon$ (black stars in Figure \ref{fractal_uncertainty_discrete}).

\subsection{Numerical computation}

For efficient numerical computations of the phase sensitivity coefficient, it is necessary to reduce the number of evaluations of the phase function, which can be computationally expensive. Toward this aim, we propose the following guidelines.

\paragraph*{Choice of the set $\mathcal{A}$.} In order to capture the overall phase sensitivity of the system, it seems natural to consider a set $\mathcal{A}$ that contains a large region of the state space. However, the phase sensitivity coefficient is unique and does not depend on the size and location of the set $\mathcal{A}$, provided that this set has a non-empty intersection with the phaseless set $\mathcal{S}$, \red{as required by our definition.} This implies that accurate results are obtained with small sets $\mathcal{A}$ that cover a small portion of $\mathcal{S}$. From a practical point of view, regions of interest containing $\mathcal{S}$ can be located through the use of finite-time Lyapunov exponents (see Figure \ref{Lyap_expo}) or simply by detecting regions of high variation of the phase function (see Figure \ref{Lorenz_global}). In addition, in the generic case where the phaseless set is of co-dimension $1$ (or less), it is more efficient to choose a set $\mathcal{A}$ of dimension $1$. Such a set has an intersection with $\mathcal{S}$ so that $\mathcal{M}_\mathcal{S}(\epsilon) \neq \emptyset$ for all $\epsilon$. It follows that the results of Section \ref{sub_sec_fractal_dim} are still valid (with the measure $\mu$ being the one-dimensional Lebesgue measure).

\paragraph*{Approximation of the ball $B$.} It is efficient to approximate the ball $B$ by considering a few sample points on its boundary. \red{In particular, two points on opposite sides of $B$ are enough to obtain the exact value of the phase sensitivity coefficient. In this case, the phase sensitivity function is approximated by
\begin{equation*}
\tilde{f}(\ve{x},\epsilon) = \max_{\ve{x}'\in \{\ve{x}_k-\epsilon \ve{e},\ve{x}_k+\epsilon \ve{e}\}} d(\Theta(\ve{x}_k),\Theta(\ve{x}')) \leq f(\ve{x},\epsilon)
\end{equation*}
where $\ve{e}$ is a unit direction. For $\ve{x} \in \mathcal{M}_\mathcal{S}(\epsilon)$,  it is clear that we have $\tilde{f}(\ve{x},\epsilon) \leq \pi$ instead of \eqref{integral1}. However, in the generic case where the phaseless set is of co-dimension $1$ (or less), the two points $\ve{x}_k-\epsilon \ve{e}$ and $\ve{x}_k+\epsilon \ve{e}$ on the boundary of $B(\ve{x},\epsilon)$ lie on both sides of $\mathcal{S}$, so that
\begin{equation}
\label{integral1b}
\lim_{\epsilon \rightarrow 0} \tilde{f}(\ve{x},\epsilon) \neq 0 \qquad \forall \ve{x} \in \mathcal{M}_\mathcal{S}(\epsilon)\,.
\end{equation}
In addition, we have
\begin{equation}
\label{integral2b}
\tilde{f}(\ve{x},\epsilon) = \epsilon |\nabla \Theta(\ve{x})\cdot \ve{e}| + \mathcal{O}(\epsilon^2) \quad \forall \ve{x} \in \mathcal{A} \cap \mathcal{B} \setminus \mathcal{M_S}(\epsilon+ \delta \epsilon)
\end{equation}
instead of \eqref{integral2}. Using \eqref{integral1b} and \eqref{integral2b} in the computations of Section \ref{sub_sec_fractal_dim}, it is easy to see that the result still holds when $f$ is replaced by $\tilde{f}$. This has been confirmed by numerical simulations.} Note also that the value of the phase sensitivity coefficient does not depend on the location of the two opposite points on the ball.

\paragraph*{\red{Method based on the sets $\mathcal{M}_{\delta \theta}$.}} As shown in Figure \ref{fractal_dim_comparison}, the phase sensitivity coefficient can also be computed with the ratio \eqref{sensitivity}. However, the value $\mu [\mathcal{M}_{\delta \theta}(\epsilon)]$ in \eqref{sensitivity} is usually underestimated when the ball $B$ is approximated with a few points (see above) so that the numerical computation may yield a zero value for small (but positive) values $\epsilon$. The computation of the phase sensitivity coefficient is therefore easier and more accurate when using the ratio \eqref{phase_uncert_coeff}.

\section{Application to bursting neurons}
\label{sec_bursting}

The notions of phase and isochrons play a central role in the study of neuron models, showing the sensitivity of neurons to external inputs. Motivated by preliminary observations presented in \cite{MRMM_bursting}, we illustrate the framework based on the phase sensitivity coefficient on popular (periodic) bursting neuron models (see Appendix \ref{app_burst_model}). Our comparison of different models \red{shows that some} elliptic bursting models exhibit strong fractal properties associated with very high phase sensitivity.

Bursting is the alternation between a relatively quiescent state and a succession of rapid spikes in a system. This phenomenon is observed with slow-fast dynamics and is explained by particular bifurcations in the fast subsystems \cite{Rinzel_Lee}. According to the types of bifurcations involved in the bursting mechanism, bursting models can be classified in different categories: elliptic (E) bursting (subcritical Hopf bifurcation/fold limit cycle bifurcation), square-wave (SW) bursting (saddle-node bifurcation/homoclinic bifurcation), parabolic (P) bursting (saddle-node on a limit cycle bifurcation); see \cite{Izhikevich_burst_classif} for more details.

\begin{figure}[h]
\begin{center}
\includegraphics[width=9cm]{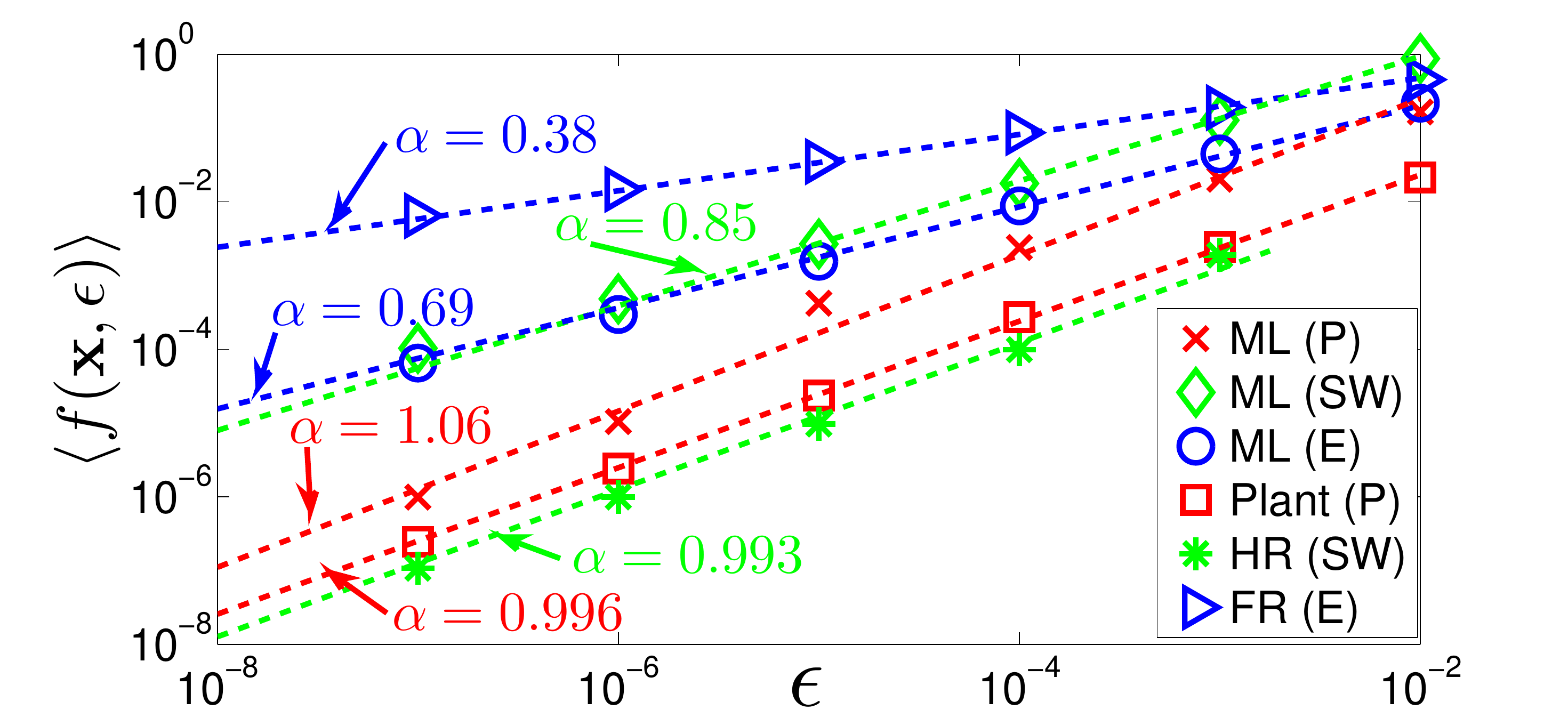}
\caption{Elliptic (E) bursting models \red{considered in our analysis} (blue) are characterized by a high phase sensitivity coefficient $\beta=1-\alpha$, which corresponds to isochrons with high fractal dimension. \red{The numbers in the figure are the slopes $\alpha$ of the curves. We obtain the following phase sensitivity coefficients $\beta$: $0.62$ (FR (E)); $0.31$ (ML (E)); $0.15$ (ML (SW)); $0.01$ (HR (SW)); $0.01$ (ML (P)); $0.00$ (Plant (P)).} Simulation details are given in Appendix \ref{app_numerical}.}
\label{fractal_bursting}
\end{center}
\end{figure}

As shown in Figure \ref{fractal_bursting}, each type of bursting is associated with a particular range of values for the phase sensitivity coefficient, an observation which reflects different values of fractal dimension and overall phase sensitivity. Elliptic (E) bursting models \red{considered in our analysis} (in blue) are characterized by a high phase sensitivity coefficient $\beta=1-\alpha$, and therefore by strong fractal properties. This is in agreement with the preliminary observations of \cite{MRMM_bursting}. In contrast, parabolic (P) bursting models (in red) have a phase sensitivity coefficient equal to $0$, thereby exhibiting no fractal properties. Square-wave (SW) bursting models (in green) roughly correspond to an intermediate situation with a low phase sensitivity coefficient. They are characterized by no fractal properties (HR model) or weak fractal properties (ML model).

\begin{figure}[h]
\begin{center}
\subfigure[]{\includegraphics[height=4.4cm]{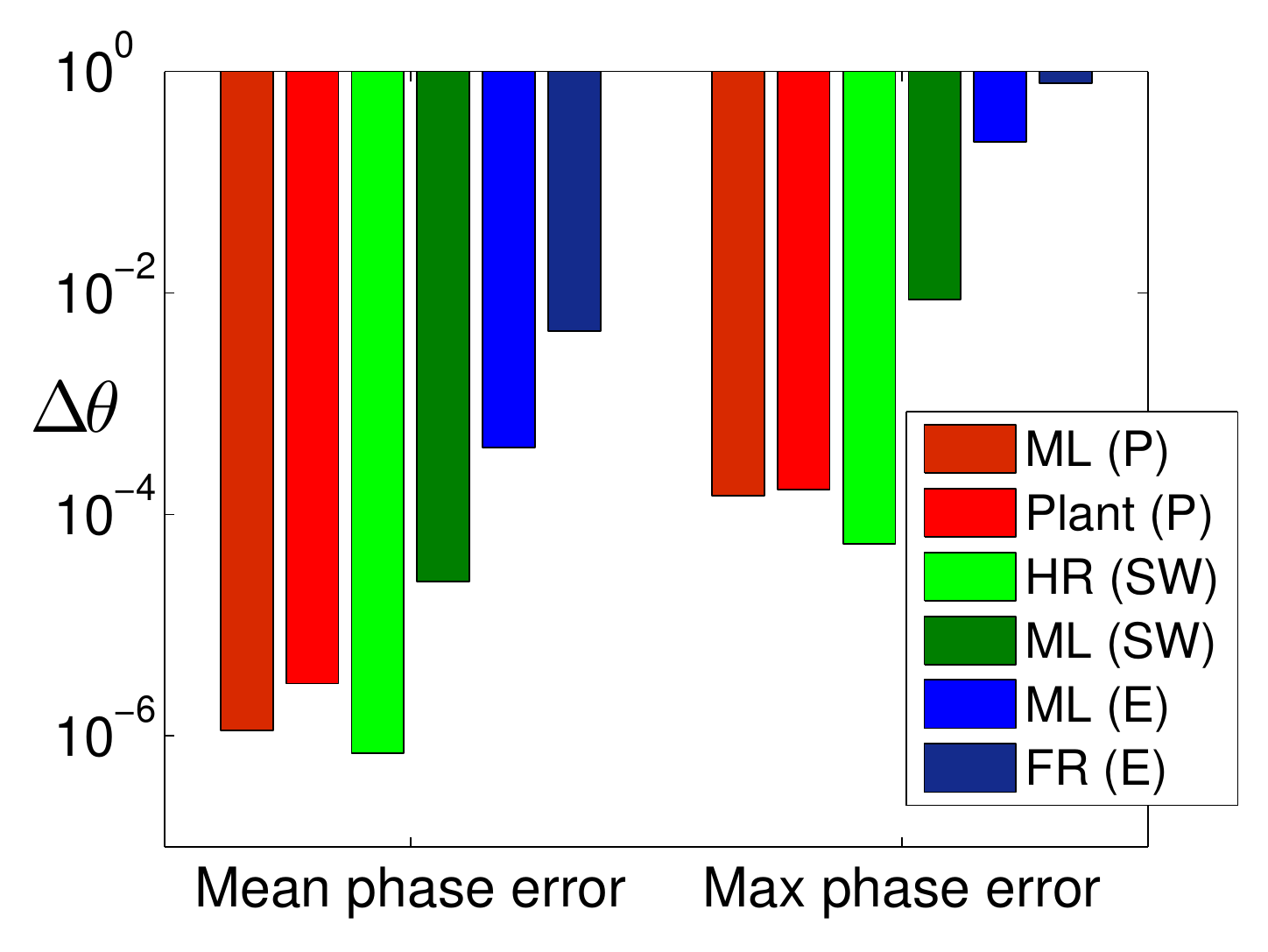}}
\subfigure[]{\includegraphics[height=4.4cm]{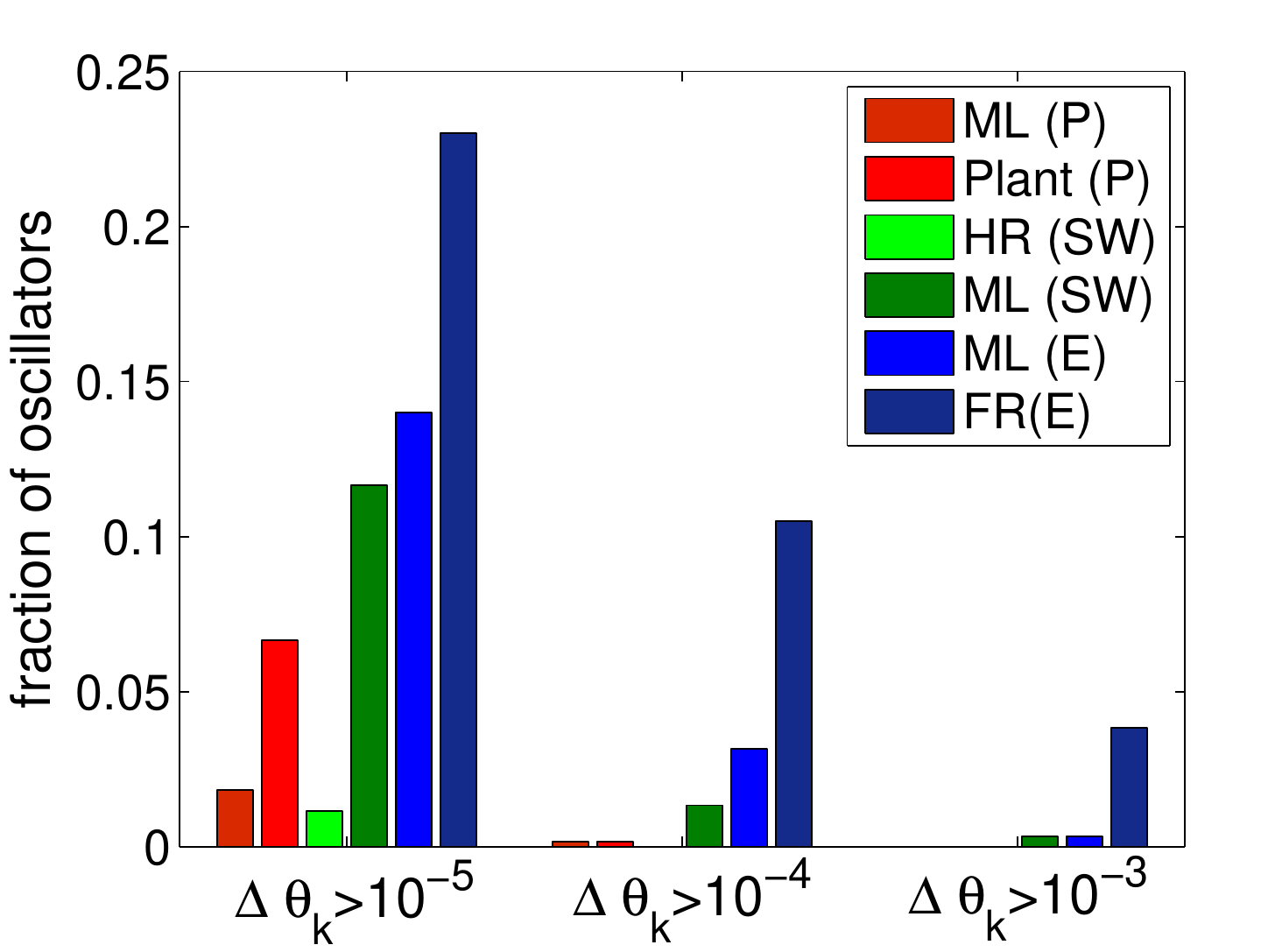}}
\caption{Networks of bursting neurons characterized by a high phase sensitivity coefficient are more sensitive to a slight perturbation of the input. The models are sorted by increasing values of the phase sensitivity coefficient. (a) Average and maximum values of the phase errors $\Delta \theta_k$ (computed over $100$ neurons and over the different simulations). \red{See also detailed results in Appendix \ref{app_network}.} (b) Fraction of oscillators whose phase error $\Delta \theta_k$ is larger than $10^{-5}$, $10^{-4}$, and $10^{-3}$. [Simulations are performed for $6$ different pulse sizes $\{0.01,0.05,0.1,0.15,0.2,0.5\} \times V_{range}$, where $V_{range}$ is the length of the interval spanned by the first state variable on the limit cycle. The perturbation $\xi_k$ follows a normal distribution of mean $0$ and standard deviation $10^{-6} \times V_{range}$.]}
\label{phase_response_network}
\end{center}
\end{figure}

Our analysis based on the phase sensitivity coefficient is confirmed by the following numerical experiment. We consider a network of $100$ neurons with random initial conditions $\ve{x}_k$ on the limit cycle (uniform distribution in phase) and \red{we assume that these neurons receive a common impulsive input $\ve{u}=\ve{e} \ \delta(t)$, where $\ve{e}$ is a vector in the $V$ direction (i.e. membrane voltage). The state of the neurons instantaneously jumps to $\ve{x}_k+\ve{e}$, which corresponds to an asymptotic phase $\theta_k$ that we compute. Then we perform the same simulation with identical initial conditions but with a slightly perturbed input $\tilde{\ve{u}}_k=(\ve{e}+\xi_k \, \ve{e}/\|\ve{e}\|) \, \delta(t)$ for each neuron, where $\xi_k \ll 1$ is a small random variable. In this case, the neurons jump to the state $\ve{x}_k+\ve{e}(1+\xi_k/\|\ve{e}\|)$ and we compute their corresponding phase $\tilde{\theta}_k$. Depending on the pulse size $\|\ve{e}\|$, the neurons may (or may not) reach a region of high phase sensitivity associated with a high phase error $\Delta \theta_k=\theta_k-\tilde{\theta}_k$. In order to cover a large part of the state space, we consider different pulse sizes and obtain statistical results on $\Delta \theta_k$ that are consistent with the values of the phase sensitivity coefficient (Figure \ref{phase_response_network}, \red{see also Appendix \ref{app_network} for detailed results obtained with different pulse sizes}).} It is also noticeable that the phase sensitivity coefficient, computed on a very small subset of the state space (Figure \ref{fractal_bursting}), captures well the sensitivity of the network, which is computed globally in a large region of the state space (Figure \ref{phase_response_network}).

\begin{figure}[h]
\begin{center}
\subfigure[]{\includegraphics[height=3.2cm]{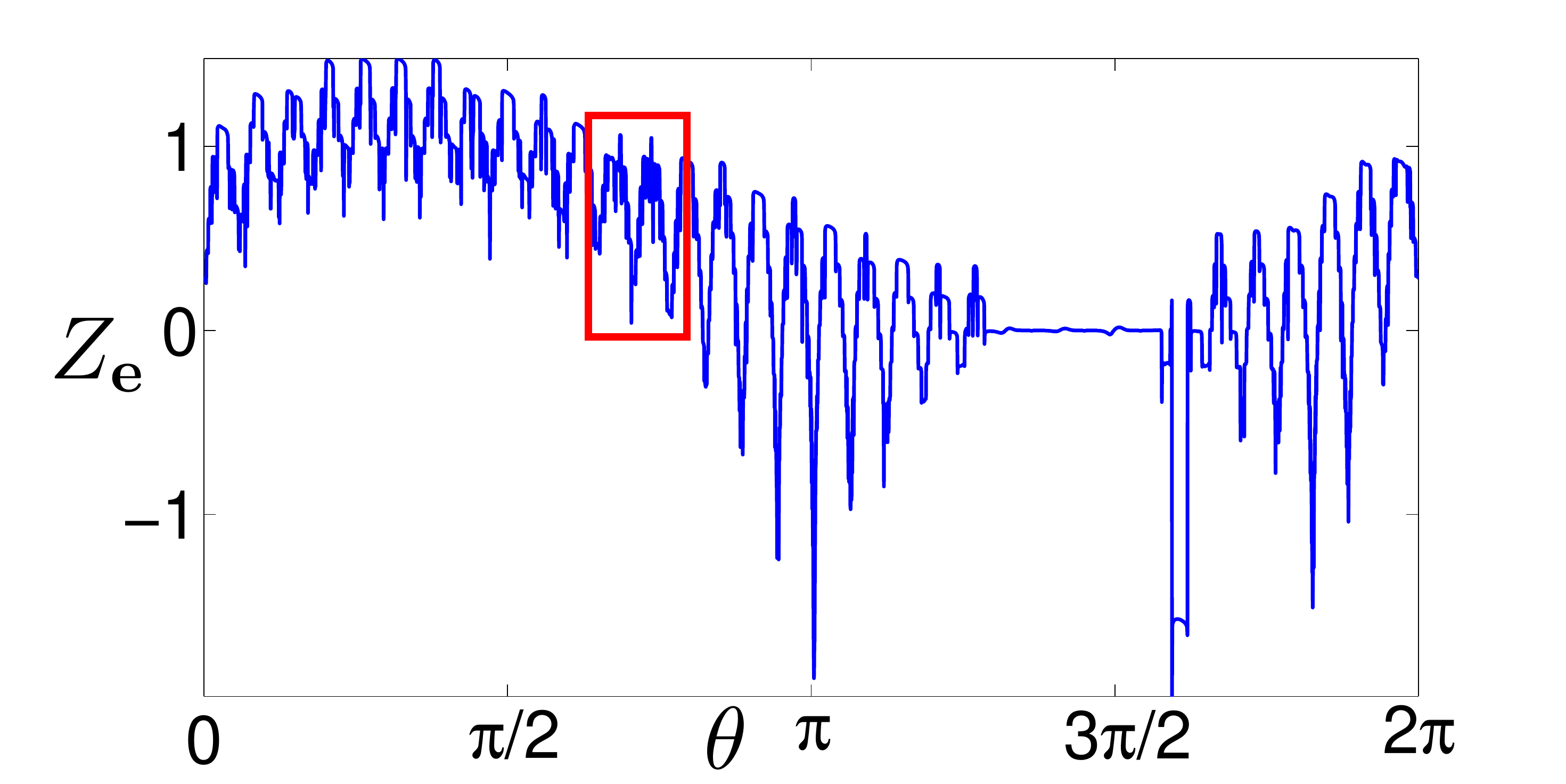}}
\subfigure[]{\includegraphics[height=3.2cm]{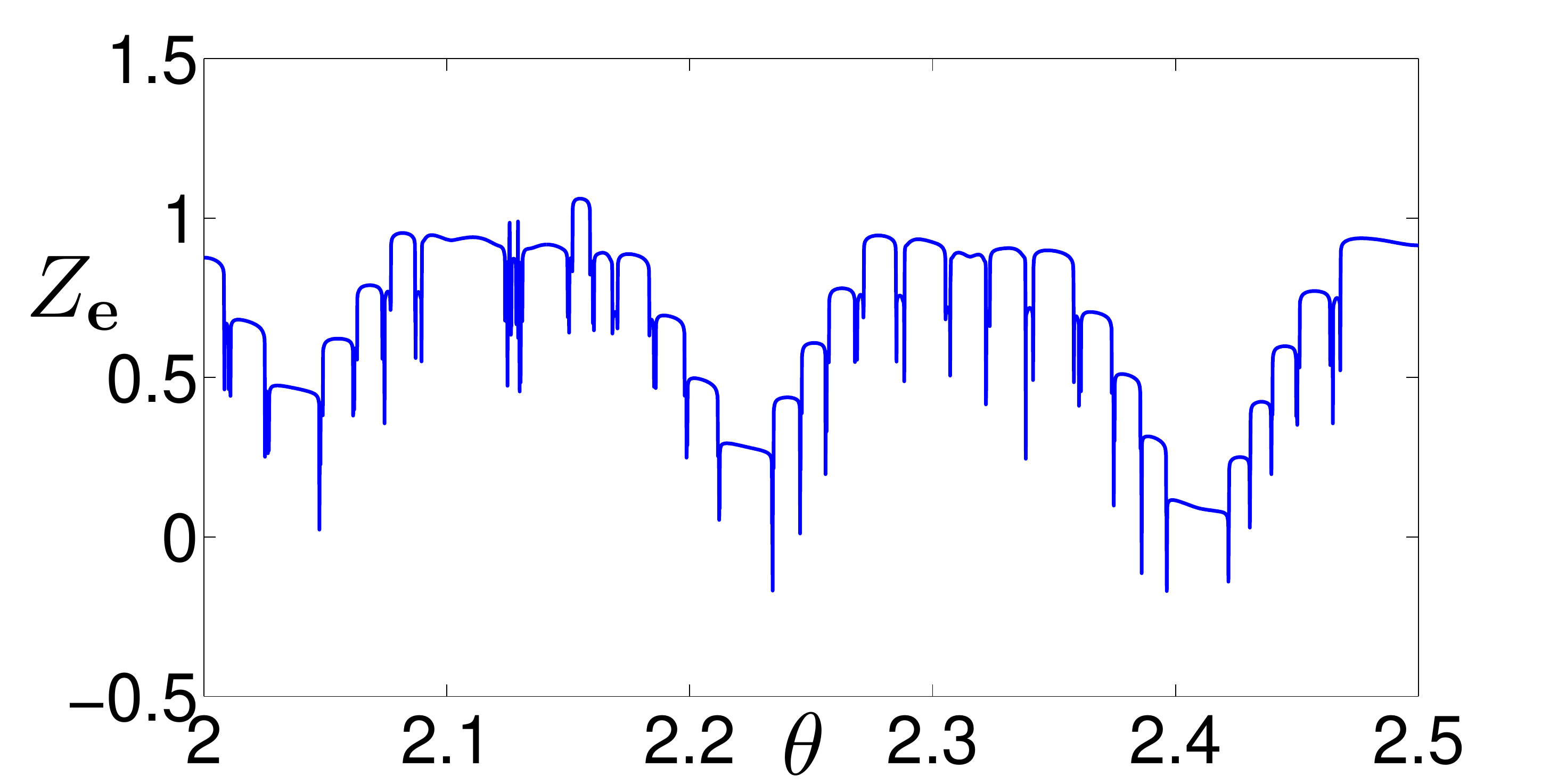}}
\caption{For the modified Morris-Lecar model in elliptic bursting regime, the phase response curve \eqref{PRC} is fractal, with a dimension equal to $1+\beta$. No fractal pattern is observed in the flat region $\theta\in[4,5]$, which corresponds to the quiescent segment on the limit cycle. The perturbation $\ve{e}$ is in the $V$ direction and has a strength of $1$ mV.}
\label{fractal_PRC_elliptic}
\end{center}
\end{figure}

Although applied to a small collection of models, the results tend to show that elliptic bursting models are characterized by a very high overall phase sensitivity. \red{For the elliptic bursting neurons tested,} a small uncertainty on the input signal (or a small noise perturbation) may induce high variation in phase (see Figure \ref{phase_response_network}), and a reduction of this uncertainty only slightly reduces the uncertainty on the phase. These neuron models are therefore characterized by sensitive and unreliable responses to inputs. \red{In particular, it is clear from the results of Figure \ref{phase_response_network} that these neurons with identical initial phases would not remain synchronized under the effect of a common (impulsive) input with a (very small) additive noise. Also, the phase response curve of these models has strong fractal properties (Figure \ref{fractal_PRC_elliptic}(a)) and cannot be used in practice since it is impossible to know exactly its values (in the fractal regions, see Figure \ref{fractal_PRC_elliptic}(b)).}
Note that the sensitivity observed here only results from the properties of the neuron dynamics and does not depend on the type of forcing (or the type of coupling in a network). This is in contrast to high sensitivity to initial conditions usually reported in neuron models with an external input (e.g. shear-induced chaos in the periodically kicked Morris-Lecar neuron \cite{Lin_shear_chaos}).

\section{Conclusion}
\label{conclu}

We have discussed the fractal properties of asymptotically periodic systems and we have provided a theoretical framework to quantify these properties. Through the notion of phase sensitivity coefficient, the fractal capacity dimension of the isochrons has been related to the overall phase sensitivity of the system.

The main implication of the results is that there exist systems---with isochrons of high fractal dimension---that are characterized by a (extremely) high phase sensitivity. For these systems, reducing the intensity of a noise perturbation only slightly decreases the average uncertainty on the phase. This is for instance the case \red{for some} elliptic bursting neuron models, whose response to external inputs is unreliable.

We will finally note that the rich geometric properties of asymptotically periodic systems illustrate the importance of developing efficient methods for computing the isochrons and the phase function. It is only when these methods are well-suited to complex dynamics in high dimensional spaces that they may unveil new properties, such as the fractal isochrons described in this paper.

\section*{Acknowledgments}
This work was funded by Army Research Office Grant W911NF-11-1-0511, with Program Manager Dr. Sam Stanton. It was completed while A. Mauroy was with the Department of Mechanical Engineering, University of California Santa Barbara. A. Mauroy is currently supported by a BELSPO (Belgian Science Policy) return grant.

\appendix

\section{Scaling of $\|\nabla \Theta\|$}
\label{app_scaling_gradient}

The (infinitesimal) phase difference between two trajectories
\begin{equation*}
\Theta(\varphi(-t,\ve{x}+d\ve{x}))-\Theta(\varphi(-t,\ve{x})) = \nabla \Theta^\textrm{T}(\varphi(-t,\ve{x})) \, (\varphi(-t,\ve{x}+d\ve{x})-\varphi(-t,\ve{x}))
\end{equation*}
is constant. It follows that we have, for all infinitesimal $d\ve{x}$,
\begin{equation*}
\nabla \Theta^\textrm{T}(\ve{x})  d\ve{x} = \nabla \Theta^\textrm{T}(\varphi(-t,\ve{x})) \, M(t) d\ve{x}\,,
\end{equation*}
where $M(t)$ is the fundamental matrix solution of
\begin{equation*}
\frac{d M}{dt}=J(\varphi(-t,\ve{x})) M
\end{equation*}
with $M(0)=I$ and $J$ is the Jacobian matrix of the vector field $\ve{F}$.

If we assume that $\mathcal{S}$ is normally hyperbolic and of co-dimension (at most) $1$, we can choose $d\ve{x}=d\ve{x}_\perp$ in the direction tangent to the fibers of $\mathcal{S}$ and we obtain
\begin{equation}
\label{gradient_theta}
\| \nabla \Theta(\ve{x}) \| \, \|d\ve{x}_\perp \| \, \cos \beta(0) = \| \nabla \Theta(\varphi(-t,\ve{x}))\| \, \|M(t) d\ve{x}_\perp \| \, \cos \beta (t)
\end{equation}
with $\beta(t)$ the angle between the gradient $\nabla \Theta(\varphi(-t,\ve{x}))$ and the direction $M(t) d\ve{x}_\perp$ tangent to the fiber. Since $\nabla \Theta(\varphi(-t,\ve{x}))\neq 0$ and $M(t) d\ve{x}_\perp \neq 0$ for all $t$, it follows that either $\cos \beta(t)=0$ for all $t$ or $\cos \beta(t) \neq 0$ for all $t$. The first case (i.e. $\beta(t)=\pm \pi/2$) is not possible since the isochrons are not the fibers of $\mathcal{S}$ (otherwise, $\mathcal{S}$ would not be phaseless). For the same reason, the isochrons are not tangent to the normal bundle of $\mathcal{S}$ in the neighborhood of $\mathcal{S}$, so that one cannot have $\cos \beta(t) \rightarrow 0$ as $t\rightarrow \infty$.  Also, one has
\begin{equation*}
\|M(t) d\ve{x}_\perp \| = \mathcal{O} \left(e^{\lambda_\perp t} \right)  \qquad t\rightarrow \infty
\end{equation*}
where $\lambda_\perp$ is the (negative) Lyapunov exponent of the system in backward time, associated with the normal direction $d\ve{x}_\perp$. (Note that the value $\exp(\lambda_\perp)$ is the generalized Lyapunov type number along the fiber \cite{Fenichel1,Wiggins}.) It follows from \eqref{gradient_theta} that
\begin{equation*}
\| \nabla \Theta(\varphi(-t,\ve{x}))\| = \mathcal{O} \left(e^{-\lambda_\perp t} \right)  \qquad t\rightarrow \infty \,.
\end{equation*}
Since the distance $d_{\mathcal{S}}$ between the trajectories and $\mathcal{S}$ satisfies $d_\mathcal{S}(\varphi(-t,\ve{x})) =\mathcal{O} (\exp(\lambda_\perp t))$ as $t \rightarrow \infty$, we have finally
\begin{equation*}
\| \nabla \Theta(\varphi(-t,\ve{x}))\| = \mathcal{O} \left(\frac{1}{d_\mathcal{S}(\varphi(-t,\ve{x}))} \right)\,.
\end{equation*}

\section{Bursting neuron models}
\label{app_burst_model}

\subsection{Modified Morris-Lecar (ML) model (square-wave, elliptic, or parabolic bursting)}

\begin{eqnarray*}
C\,\dot{V} & = & -g_{Ca} m_\infty(V)(V-V_{Ca})-  g_K n (V-V_K) - g_L (V-V_L) \\
 & &- g_{KCa} z(h) (V-V_K) - g_{CaS} s (V-V_{Ca}) + I \\
\dot{n} & = & \phi\, (w_\infty(V)-n)/\tau(V) \\
\dot{h} & = & \epsilon_1(-\mu g_{Ca} m_\infty(V) \, (V-V_{Ca}) - h) \\
\dot{s} & = & \epsilon_2(s_\infty(V)-s)/\tau_s \\
\end{eqnarray*}
with
\begin{eqnarray*}
m_\infty(V) & = & 0.5\left(1+\tanh\frac{V-V_1}{V_2}\right)  \\
w_\infty(V) & = & 0.5\left(1+\tanh\frac{V-V_3}{V_4}\right)  \\
s_\infty(V) & = & 0.5\left(1+\tanh\frac{V-V_5}{V_6}\right) \\
z(h) & = & \frac{h}{\textrm{Ca}_0+h} \,, \\
\tau(V) & = & \cosh^{-1}\frac{V-V_3}{2 V_4}\,.
\end{eqnarray*}
The parameters are given in Table \ref{table1}.

\begin{table}
	\centering
		\begin{tabular}{lccc}
			Parameter & Square-wave & Elliptic & Parabolic \\
			\hline
			$g_{Ca}$ & 4 & 4 & 4 \\
$g_K$ ($\textrm{mS/cm}^2$)  & 8 & 8 & 8 \\
$g_L$ ($\textrm{mS/cm}^2$) & 2 & 2 & 2 \\
$V_K$ (mV) & -84 & -84 & -84 \\
$V_L$ (mV)& -60 & -60 & -60 \\
$V_{Ca}$ (mV) & 120 & 120 & 120 \\
$C$ ($\mu\textrm{F/cm}^2$) & 17.8 & 10 & 1 \\
$I$ ($\mu \textrm{A/cm}^2$) & 45 & 120 & 65 \\
$g_{KCa}$ ($\textrm{mS/cm}^2$) & 0.25 & 0.75 & 1 \\
$\phi$ & 0.25 & 0.04 & 1.333 \\
$\epsilon_1$ & 0.005 & 0.002 & 0.02 \\
$\mu$ & 0.2 & 0.3 & 0.025 \\
$\epsilon_2$ & 0 & 0 & 0.02 \\
$\tau_s$ (ms) & - & - & 0.05 \\
$g_{Cas}$ ($\textrm{mS/cm}^2$) & - & - & 1 \\
$\textrm{Ca}_0$ & 10 & 18 & 1 \\
$V_1$ (mV) & -1.2 & -1.2 & -1.2 \\
$V_2$ (mV) & 18 & 18 & 18 \\
$V_3$ (mV) & 12 & 2 & 12 \\
$V_4$ (mV) & 17.4 & 30 & 17.4 \\
$V_5$ (mV) & - & - & 12 \\
$V_6$ (mV) & - & - & 24
		\end{tabular}
		\caption{Parameters of the modified Morris-Lecar (ML) model}
		\label{table1}
\end{table}

\subsection{Hindmarsh-Rose (HR) model (square-wave bursting)}

\begin{eqnarray*}
\dot{V} & = & n-aV^3+bV^2-h+I \\
\dot{n} & = & c-dV^2-n  \\
\dot{h} & = & r(\sigma(V-V_0)-h) 
\end{eqnarray*}
The parameters are $a=1$, $b=3$, $c=1$, $d=5$, $r=0.001$, $\sigma=4$, $V_0=-1.6$, $I=2$.

\subsection{FitzHugh-Rinzel (FR) model (elliptic bursting)}

\begin{eqnarray*}
\dot{V} & = & V-V^3/3-w+y+I \\
\dot{w} & = & \delta(a+V-bw) \\
\dot{y} & = & \mu(c-V-dy)
\end{eqnarray*}
The parameters are $I = 0.3125$, $a = 0.7$, $b = 0.8$, $c = -0.9$, $d = 1$, $\delta = 0.08$, and $\mu =0.001$.

\subsection{Plant model (parabolic bursting)}

\begin{eqnarray*}
C\,\dot{V} & = & -g_{Na} m^3_\infty(V) h (V-V_{Na}) -  g_{Ca} x (V-V_{Ca}) \\
 & & -\left(g_K n^4 +\frac{k_{KCa}\, c}{0.5+c} \right) (V-V_K) - g_L (V-V_L) \\
\dot{n} & = & (h_\infty(V)-h)/\tau_h(V) \\
\dot{h} & = & (n_\infty(V)-n)/\tau_n(V) \\
\dot{x} & = & (x_\infty(V)-x)/\tau_x \\
\dot{c} & = & f \, (k_1 x (V_{Ca}-V)-c)
\end{eqnarray*}
with
\begin{eqnarray*}
w_\infty(V) & = & \frac{\alpha_\infty(V)}{\alpha_\infty(V)+\beta_\infty(V)} \quad \textrm{for } w=m,h,n\\
\tau_\infty(V) & = & \frac{12.5}{\alpha_\infty(V)+\beta_\infty(V)} \quad \textrm{for } w=h,n\\
x_\infty(V) & = & \frac{1}{\exp(-0.15(V+50))+1}
\end{eqnarray*}
where
\begin{eqnarray*}
\alpha_m(V) & = & 0.1\frac{50-V_s}{\exp((50-V_s)/10)-1} \\
\beta_m(V) & = & 4 \exp((25-V_s)/18) \\
\alpha_n(V) & = & 0.01 \frac{55-V_s}{\exp((55-V_s)/10)-1} \\
\beta_n(V) & = & 0.125 \exp((45-V_s)/80) \\
\alpha_h(V) & = & 0.07 \exp((25-V_s)/20) \\
\beta_h(V) & = & \frac{1}{\exp((55-V_s)/10)+1}
\end{eqnarray*}
and $V_S=127/105 V+8265/105$. The parameters are $C=1\, \mu\textrm{F/cm}^2$, $g_{Ca}=0.004\, \textrm{mS/cm}^2$, $g_{Na}=4\, \textrm{mS/cm}^2$, $g_{K}=0.3\, \textrm{mS/cm}^2$, $g_{L}=0.004\, \textrm{mS/cm}^2$, $f=0.0003\, \textrm{ms}^{-1}$, $g_{KCa}=0.03\, \textrm{mS/cm}^2$, $V_{Ca}=140\, \textrm{mV}$, $V_{Na}=30\, \textrm{mV}$, $V_{K}=-75\, \textrm{mV}$, $V_{Ca}=-40\, \textrm{mV}$, and $k_1=0.0085\, \textrm{mV}^{-1}$.


\section{Numerical details and simulation parameters}
\label{app_numerical}

\subsection{Computation of the phase function}

We computed the Fourier averages \eqref{Fourier_av}-\eqref{Fourier_av2} over a finite time horizon $T$. For continuous-time models, the integral was obtained with the MATLAB function \guillemets{trapz} and the trajectories were computed with the MATLAB function \guillemets{ode45}. For the bursting neuron models (Section \ref{sec_bursting}), we integrated over only one limit cycle period $T_0=2\pi/\omega_0$, considering the truncated Fourier averages
\begin{equation*}
\Theta(\ve{x}) \approx \angle \left(\frac{1}{T_0} \int_{T-T_0}^T g \circ \varphi(t,\ve{x}) \, e^{-i\omega_0 t} \, dt \right) \,.
\end{equation*}
The simulations parameters are summarized in Table \ref{table_simu}.

\begin{table}[h!]

\centering
		\begin{tabular}{cccccc}
		Model & Frequency  & Relative error  & Absolute error & Time horizon  & Function $g$ \\
		& $\omega_0$ & (for ode45) & (for ode45) & $T$ & \\
		\hline
		Van der Pol & 0.942958 & 1e-6 & 1e-50 & 100 & $y$ \\
		Lorenz & 15.4547 & 1e-9 & 1e-300 & 50 & $x$ \\
		ML square-wave & 0.008870246 & 1e-6 & 1e-300 & 3500 & $n$ \\
		ML elliptic & 0.0037015 & 1e-6 & 1e-300 & 8500 & $h$ \\
		ML parabolic & 0.075131 & 1e-6 & 1e-300 & 500 & $n$ \\
		HR & 0.014586 & 1e-6 & 1e-300 & 2000 & $n$ \\
		FR & 0.008218 & 1e-6 & 1e-300 & 4000 & $y$ \\
		Plant & 0.00058225 & 1e-6 & 1e-300 & 30000 & $h$ \\
		Discrete-time map \eqref{discrete_map} & 1.53828241 & - & - & 5000 & y \\
		Discrete-time map of Figure \ref{vdp_iso}(b) & 3.52690624 & - & - & 5000 & y
					\end{tabular}
					\caption{Simulation parameters for the computation of the phase function}
					\label{table_simu}
\end{table}

\subsection{Computation of the phase sensitivity coefficient}

The phase sensitivity function was computed with two points on the boundary of the balls $B(\ve{x},\epsilon)$ and its average was obtained with a finite number $n_{pt}$ of sample points equally distributed in $\ve{x}_k \in \mathcal{A}$. We therefore considered the approximation
\begin{equation*}
\langle f(\ve{x},\epsilon) \rangle_{\mathcal{A} \cap \mathcal{B}} \approx \frac{1}{n_{pt}} \sum_{k=1}^{n_{pt}} \max_{\ve{x}'\in \{\ve{x}_k-\epsilon \ve{e},\ve{x}_k+\epsilon \ve{e}\}} d(\Theta(\ve{x}_k),\Theta(\ve{x}'))
\end{equation*}
where $\ve{e}$ is a unit vector.

Figures \ref{phase_uncertainty}, \ref {fractal_dim_comparison}, and \ref{fractal_bursting} were obtained with the parameters given in Table \ref{table_simu2}.

\begin{table}[h!]
\centering

		\begin{tabular}{cccc}
		Model & Set $\mathcal{A}$  & Number of sample points $n_{pt}$ & direction of $\ve{e}$\\
		\hline
		Van der Pol & $[-0.5,0.5] \times \{0\}$ & 1000 & along $x$\\
		Lorenz & $[-48.8,-48.75] \times \{100\} \times \{319\}$ & 2500 & along $x$\\
		ML square-wave & $\{-15\} \times [0.1,0.2] \times \{12\}$ & 10000 & along $V$ \\
		ML elliptic & $\{30\} \times [0,0.5]  \times \{16\}$ & 10000 & along $V$\\
		ML parabolic & $\{0\} \times [0.2,0.3] \times \{1.5\} \times \{0.15\}$ & 10000 & along $V$\\
		HR & $\{0.5\} \times [-10,4] \times \{1.9\}$ & 2500 & along $V$\\
		FR & $\{-1\} \times [-0.5,0.5] \times \{0.01\}$ & 10000 & along $V$\\
		Plant & $\{-20\} \times [0,1] \times \{0.4\} \times \{0.74\} \times \{0.6\}$ & 10000 & along $V$\\
		Discrete-time maps & $y\in[0,1]$ & 10000 & along $y$
					\end{tabular}
					\caption{Simulation parameters for the computation of the phase sensitivity coefficient}
					\label{table_simu2}
\end{table}

\red{
\section{Detailed results related to Figure \ref{phase_response_network}}
\label{app_network}
}

\red{The following tables contain detailed results related to the numerical experiment summarized in Figure \ref{phase_response_network}(a) (Section \ref{sec_bursting}). For different pulse sizes $\|\ve{e}\|$, Table \ref{table_network_mean} and Table \ref{table_network_max} show the mean phase error and the maximum phase error, respectively.
}

\begin{table}[h!]
\red{
\centering
		\begin{tabular}{ccccccccccccc}
		Pulse Size $\|\ve{e}\|$ && ML(P) && Plant (P) && HR (SW) && ML (SW) && ML (E) && FR (E) \\
		\hline
 $0.01 \times V_{range}$	  && 1.06e-06 &&  8.14e-06  && 2.42e-07 &&  3.92e-06 &&  2.37e-03 &&  \textbf{9.14e-03} \\
 $0.05 \times V_{range}$  && 5.98e-07  && 2.97e-06 &&  2.54e-07  && 6.85e-06 &&  1.91e-05 &&  \textbf{1.40e-02} \\
 $0.1 \times V_{range}$  && 1.42e-06 &&  1.05e-06  && 3.00e-07 &&  7.39e-06 &&  7.97e-07 &&  \textbf{3.83e-03} \\
 $0.15 \times V_{range}$  && 2.98e-06  && 2.16e-06  && 1.02e-06  && 4.73e-06  && 1.43e-06 &&  \textbf{8.36e-05} \\
 $0.2 \times V_{range}$  && 4.50e-07  && 3.32e-06  && 4.53e-07  && 2.04e-06  && 1.99e-06  && \textbf{1.20e-04} \\
 $0.5 \times V_{range}$  && 1.95e-07  && 2.13e-07 &&  1.97e-06  && \textbf{1.97e-05}  && 1.08e-06  && 1.41e-06 \\
\hline
Mean && 1.12e-06 &&  2.97e-06 &&  7.06e-07 &&  7.44e-06 &&  3.99e-04 &&  \textbf{4.53e-03}
					\end{tabular}
					\caption{Mean phase error. The maximum value for a given pulse size is in bold.}
					\label{table_network_mean}
					}
\end{table}

\begin{table}[h!]
\red{
\centering
		\begin{tabular}{ccccccccccccc}
		Pulse Size $\|\ve{e}\|$ && ML(P) && Plant (P) && HR (SW) && ML (SW) && ML (E) && FR (E) \\
		\hline
$0.01 \times V_{range}$ &&   3.50e-05 &&  1.67e-04 &&  8.51e-06  && 3.02e-05 &&  2.32e-01 &&  \textbf{6.48e-01} \\
 $0.05 \times V_{range}$  &&  1.43e-05 &&  6.59e-05 &&  4.56e-06  && 2.09e-04 &&  8.59e-04  &&\textbf{ 7.78e-01} \\
 $0.1 \times V_{range}$  &&  2.83e-05  && 7.62e-06  && 3.44e-06  && 2.23e-04  && 1.73e-05  && \textbf{3.79e-01} \\
 $0.15 \times V_{range}$  &&  1.48e-04  && 3.64e-05 &&  5.41e-05  && 4.49e-05  && 7.45e-05  && \textbf{7.73e-03} \\
 $0.2 \times V_{range}$  &&  1.90e-05  && 4.26e-05  && 7.23e-06  && 3.43e-05  && 1.32e-04  && \textbf{1.04e-02} \\
 $0.5 \times V_{range}$  &&  2.66e-06  && 1.54e-05  && 4.85e-05  && \textbf{1.06e-03}  && 2.99e-05  && 6.21e-05 \\
\hline
Maximum && 1.48e-04  && 1.67e-04  && 5.41e-05 &&  1.06e-03  && 2.32e-01 &&  \textbf{7.78e-01}
					\end{tabular}
					\caption{Maximum phase error. The maximum value for a given pulse size is in bold.}
					\label{table_network_max}
					}
\end{table}

\bibliographystyle{siam}

\end{document}